\documentclass[a4paper,10.5pt, reqno, oneside]{amsart}

\usepackage{amssymb}
\usepackage{amstext}
\usepackage{amsmath}
\usepackage{amscd}
\usepackage{amsthm}
\usepackage{amsfonts}
\usepackage{enumerate}
\usepackage{graphicx}
\usepackage{latexsym}
\usepackage{mathrsfs}
\usepackage{mathtools}
\usepackage{tikz}
\usetikzlibrary{arrows}
\usetikzlibrary{snakes}
\usetikzlibrary{shapes}
\usetikzlibrary{calc,decorations.markings}

\usepackage{tikz-cd}

\usepackage{hyperref}

\usepackage{pstricks}
\usepackage{lscape}
\usepackage{comment}

\makeatletter
  
  \@addtoreset{equation}{section}
\makeatother

\newtheorem{theorem}{Theorem}[section]
\newtheorem{lemma}[theorem]{Lemma}
\newtheorem{proposition}[theorem]{Proposition}
\newtheorem{corollary}[theorem]{Corollary}

\newtheorem{question}[theorem]{Question}

\theoremstyle{definition}
\newtheorem{definition}[theorem]{Definition}
\newtheorem{definiton-theorem}[theorem]{Definition-Theorem}
\newtheorem{lemma-definition}[theorem]{Lemma-Definition}
\newtheorem{example}[theorem]{Example}

\theoremstyle{remark}
\newtheorem{remark}[theorem]{Remark}

\newcommand{\rmZ}{\mathrm{Z}}

\newcommand{\fkm}{\mathfrak{m}}

\newcommand{\ZZ}{\mathbb{Z}}

\newcommand{\RR}{\mathbb{R}}

\newcommand{\sfA}{\mathsf{A}}

\newcommand{\sfM}{\mathsf{M}}
\newcommand{\sfN}{\mathsf{N}}

\newcommand{\sfR}{\mathsf{R}}
\newcommand{\sfT}{\mathsf{T}}

\newcommand{\rank}{\operatorname{rank}}

\newcommand{\Hom}{\operatorname{Hom}}

\newcommand{\SL}{\operatorname{SL}}

\newcommand{\charac}{\operatorname{char}}

\newcommand{\Cl}{\operatorname{Cl}}

\newcommand{\End}{\operatorname{End}}
\newcommand{\gldim}{\mathrm{gl.dim}}

\newcommand{\add}{\mathsf{add}}
\newcommand{\CM}{\mathsf{CM}}

\begin{document}
%\tikzset{auto}

%%%%%---title---%%%%%
\title[Remarks on divisorial ideals arising from dimer models]{Remarks on divisorial ideals arising from dimer models}
\author[Yusuke Nakajima]{Yusuke Nakajima}
\date{}
%\thanks{The author is supported by Grant-in-Aid for JSPS Fellows No. 26-422.}


\address[Yusuke Nakajima]{Graduate School Of Mathematics, Nagoya University, Furocho, Chikusaku, Nagoya, 464-8602 Japan} 
\email{m06022z@math.nagoya-u.ac.jp}

\maketitle

%%%---abstract---%%%
\begin{abstract}
The Jacobian algebra $\mathsf{A}$ arising from a consistent dimer model is derived equivalent to 
crepant resolutions of a $3$-dimensional Gorenstein toric singularity $R$, and it is also called a non-commutative crepant resolution of $R$. 
This algebra $\mathsf{A}$ is a maximal Cohen-Macaulay (= MCM) module over $R$, and it is a finite direct sum of rank one MCM $R$-modules. 
In this paper, we observe a relationship between properties of a dimer model 
and those of MCM modules appearing in the decomposition of $\mathsf{A}$ as an $R$-module. 
More precisely, we take notice of isoradial dimer models and divisorial ideals which are called conic. 
Especially, we investigate them for the case of $3$-dimensional Gorenstein toric singularities associated with reflexive polygons. 
%and obtain several examples of rank one MCM modules which are not conic. 
\end{abstract}

%\setcounter{tocdepth}{1}
%\tableofcontents

%%%%%---text_start---%%%%%
\section{Introduction}
In this paper, we discuss properties of rank one MCM modules arising from non-commutative crepant resolutions of a $3$-dimensional Gorenstein toric singularity. 
Firstly, we recall the definition of non-commutative crepant resolutions \cite{VdB2}. 

\begin{definition}
Let $S$ be a $d$-dimensional local CM normal domain, and $N$ be a non-zero reflexive $S$-module. 
We say $\Lambda\coloneqq\End_S(N)$ is a \emph{non-commutative crepant resolution (= NCCR)} of $S$ if $\gldim\,\Lambda=d$ and  
$\Lambda$ is a maximal Cohen-Macaulay (= MCM) $S$-module. 
\end{definition}

We can obtain NCCRs of a $3$-dimensional complete local Gorenstein toric singularity via ``consistent dimer models". 
A dimer model $\Gamma$ is a polygonal cell decomposition of the two-torus whose vertices and edges form a finite bipartite graph, 
and we can obtain the quiver $Q_{\Gamma}$ as the dual of a dimer model $\Gamma$. 
From this quiver $Q_{\Gamma}$, we define the complete Jacobian algebra $\sfA_{Q_\Gamma}$. 
Under a certain condition so-called consistency condition, $\sfA_{Q_\Gamma}$ is an NCCR of the center $\rmZ(\sfA_{Q_\Gamma})$ and 
such a center is a $3$-dimensional complete local Gorenstein toric singularity. 
On the other hand, for each $3$-dimensional complete local Gorenstein toric singularity, there is a consistent dimer model such that
the center of the complete Jacobian algebra is isomorphic to a given toric singularity. 
Therefore, an NCCR always exists. Also, it is known that NCCRs are derived equivalent to the ordinary crepant resolutions. 
(For more details, see subsection~\ref{subsec_dimer}.) 

\medskip

In what follows, we suppose $R$ is a $3$-dimensional complete local Gorenstein toric singularity. 
Note that the Krull-Schmidt condition holds for $R$ in our situation. 
Let $\Gamma$ be a consistent dimer model which gives an NCCR of $R$ as the complete Jacobian algebra $\sfA_{Q_\Gamma}$. 
Since $\sfA_{Q_\Gamma}$ is an NCCR of $R$, it is an MCM $R$-module. 
Thus, it is decomposed as the direct sum of MCM $R$-modules: 
\begin{eqnarray}
\label{decomp_NCCR}
\sfA_{Q_\Gamma}\cong M_0^{\oplus a_0}\oplus M_1^{\oplus a_1}\oplus\cdots\oplus M_r^{\oplus a_r}. 
\end{eqnarray}
Moreover, we have $\rank_RM_i=1$ for $i=0,1,\cdots,r$ and $R$ is contained in $\add_R\sfA_{Q_\Gamma}$ 
by a construction of an NCCR (see Theorem~\ref{NCCR1}), thus let $M_0\coloneqq R$. 

In this way, we obtain some rank one MCM $R$-modules from a consistent dimer model. 
Rank one MCM modules over toric singularities are investigated in several papers e.g. \cite{BG1,Per1,Per2}. 
Especially, the number of rank one MCM modules is finite up to isomorphism \cite[Corollary~5.2]{BG1}. 
With these backgrounds, we will consider the following question. 

\begin{question}
Is there a relationship between a property of consistent dimer models and that of MCM $R$-modules appearing in the decomposition $(\ref{decomp_NCCR})$? 
\end{question}

In order to consider this question, we first observe some well-known cases as in Examples~\ref{ex_conifold} and \ref{ex_simplicial}. 
In these examples, divisorial ideals which are called conic (see subsection~\ref{toric_pre}) appear in the decomposition (\ref{decomp_NCCR}). 
As we will see later, conic divisorial ideals are always MCM modules, and related with the structure of the Frobenius push-forward of $R$ if $\charac R>0$. 
Meanwhile, dimer models appearing in such examples satisfy the isoradial (or geometrically consistency) condition.  
This condition is stronger than the consistency condition, and such a dimer model exists for every $3$-dimensional Gorenstein toric singularity by \cite{Gul}. 
Thus, we concern a relationship between conic divisorial ideals and isoradial dimer models. 
In this paper, we investigate them for the case of $3$-dimensional Gorenstein toric singularities associated with reflexive polygons 
(see subsection~\ref{subsec_reflexive}). 
For these cases, divisorial ideals arising from consistent dimer models are computed in \cite{Nak}. 
Thus, by checking divisorial ideals appearing in \cite{Nak}, we have the following results. 
(For further details on terminologies, see later sections.)

\begin{theorem}
\label{main}
Let $R$ be a $3$-dimensional complete local Gorenstein toric singularity whose associated polygon $\Delta\subset\RR^2$ (see Remark~\ref{toric_diagram}) 
is a reflexive polygon. Suppose $Q_\Gamma$ is the quiver associated with a consistent dimer model $\Gamma$ such that $R\cong\rmZ(\sfA_{Q_\Gamma})$, 
and $\{R, M_1, \cdots ,M_r\}$ are rank one MCM $R$-modules arising from the complete Jacobian algebra $\sfA_{Q_\Gamma}$ as in $(\ref{decomp_NCCR})$. 
Then we have the following. 
\begin{itemize}
\item [$(1)$] If $\Gamma$ is an isoradial dimer model, then $\{R, M_1, \cdots ,M_r\}$ are all conic divisorial ideals. 
\item [$(2)$] Conversely, if $\{R, M_1, \cdots ,M_r\}$ are all conic divisorial ideals, then $\Gamma$ is isoradial. 
\end{itemize}
\end{theorem}

By this theorem, we could obtain a relationship between conic divisorial ideals and isoradial dimer models for some cases. 
However, we remark that every conic divisorial ideal does not necessarily arise from isoradial dimer models 
(see the case of type 6a in subsection~\ref{subsec_reflexive}).

\medskip

Furthermore, we have the next corollary. 

\begin{corollary}
\label{main_cor}
Let $R$ be the same as Theorem~\ref{main}. 
If there exists a consistent dimer model associated with $R$ which is not isoradial, then there is a rank one MCM module which is not conic. 
\end{corollary}

The existence of a rank one MCM module which is not conic was investigated by Bae\c{t}ica and Bruns \cite{Bae, Bru}, 
and they gave such a module by using the Segre product. 
On the other hand, we obtain such a module by the difference between an isoradial dimer model and a non-isoradial consistent dimer model. 
Also, the existence of such a module was discussed in the context of the structure of Frobenius push-forward \cite[Question~3.2]{Wat}. 

\medskip

Thus, it is natural to ask the following, but we do not know an answer for now. 

\begin{question}
Can we obtain statements as in Theorem~\ref{main} and Corollary~\ref{main_cor} for any $3$-dimensional Gorenstein toric singularities?
\end{question}

The content of this paper is the following. 
In Section~\ref{sec_pre}, we prepare some basic results about toric singularities and dimer models to be used throughout.  
After that, in Section~\ref{sec_conic_dimer}, we observe $3$-dimensional Gorenstein toric singularities associated with reflexive polygons, 
and show Theorem~\ref{main} and Corollary~\ref{main_cor}. 

\subsection*{Notations}  
Throughout, we assume that $k$ is an algebraically closed field and $\charac k=0$ or $p\gg 0$. 
For a commutative ring $R$, we denote $\CM R$ to be the full subcategory consisting of MCM $R$-modules, 
$\add_R M$ to be the full subcategory consisting of direct summands of finite direct sums of some copies of an $R$-module $M$. 
We say $M$ is a \emph{generator} if $R\in\add_R M$. 
We denote the $R$-dual functor by $(-)^*=\Hom_R(-,R)$. 
Also, we denote by $\Cl(R)$ the class group of $R$. 
When we consider a rank one reflexive $R$-module $I$ as an element of $\Cl(R)$, we denote it by $[I]$.

%%%%%%%%%%%%%%%%%%%%%%%%%%%%%%%%%%%%%%%%%%%%%%%%%%%%%%%%%%%%
\section{Preliminaries}
\label{sec_pre}

\subsection{Preliminaries of toric singularities} 
\label{toric_pre}
In this subsection, we recall some basic facts about toric singularities and their divisorial ideals. 
For more details, see textbooks e.g. \cite{BG2, CLS}. 

Let $\sfN\cong\ZZ^d$ be a lattice of rank $d$ and let $\sfM\coloneqq\Hom_\ZZ(\sfN, \ZZ)$ be the dual lattice of $\sfN$. 
We set $\sfN_\RR\coloneqq\sfN\otimes_\ZZ\RR$ and $\sfM_\RR\coloneqq\sfM\otimes_\ZZ\RR$ and 
denote an inner product by $\langle\;,\;\rangle:\sfM_\RR\times\sfN_\RR\rightarrow\RR$. 
Let 
\[
\sigma\coloneqq\mathrm{Cone}(v_1, \cdots, v_n)=\RR_{\ge 0}v_1+\cdots +\RR_{\ge 0}v_n
\subset\sfN_\RR 
\]
be a strongly convex rational polyhedral cone of rank $d$ 
generated by $v_1, \cdots, v_n\in\ZZ^d$. We assume this system of generators is minimal. 
We consider the dual cone $\sigma^\vee$: 
\[
\sigma^\vee\coloneqq\{x\in\sfM_\RR\,|\,\langle x,y\rangle\ge0 \text{ for all } y\in\sigma \}. 
\]
Then $\sigma^\vee\cap\sfM$ is a positive affine normal semigroup. 
We define the $\fkm$-adic completion of a toric singularity 
%(or  affine normal semigroup ring) 
\[
R\coloneqq k[[\sigma^\vee\cap\sfM]]=k[[t_1^{a_1}\cdots t_d^{a_d}\mid (a_1, \cdots, a_d)\in\sigma^\vee\cap\sfM]] 
\] 
where $\fkm$ is the irrelevant maximal ideal. 
Then $R$ is a $d$-dimensional Cohen-Macaulay normal domain. 
For each generator, we define a linear form $\lambda_i(-)\coloneqq\langle-, v_i\rangle$, 
and denote $\lambda(-)\coloneqq(\lambda_1(-),\cdots,\lambda_n(-))$. 
For each $\mathbf{u}=(u_1, \cdots, u_n)\in\RR^n$, we set 
\[
\mathbb{T}(\mathbf{u})\coloneqq\{x\in\sfM \mid (\lambda_1(x), \cdots, \lambda_n(x))\ge(u_1, \cdots, u_n)\}. 
\]
Then we define the divisorial ideal $T(\mathbf{u})$ generated by all monomials whose exponent vector is in $\mathbb{T}(\mathbf{u})$. 
By the definition, we have $T(\mathbf{u})=T(\ulcorner \mathbf{u}\urcorner)$ where $\ulcorner \; \urcorner$ implies the round up 
and $\ulcorner \mathbf{u}\urcorner=(\ulcorner u_1\urcorner, \cdots, \ulcorner u_n\urcorner)$, 
hence we assume $\mathbf{u}\in\ZZ^n$ in the rest of this paper. 
%Clearly, we have $\mathbb{T}(0)=\sigma^\vee\cap\sfM$ and $T(0)=R$. 
Set $I_i\coloneqq T(\delta_{i1},\cdots,\delta_{in})$ where $\delta_{ij}$ is the Kronecker delta. 
Then a divisorial ideal $T(u_1,\cdots,u_n)$ corresponds to the element $u_1[I_1]+\cdots+u_n[I_n]$ in $\Cl(R)$. 
In general, a divisorial ideal of $R$ takes this form (see \cite[Theorem~4.54]{BG2}) and we have the following.

\begin{lemma}(see e.g. \cite[Corollary~4.56]{BG2})
\label{div_eq}
For $\mathbf{u}, \mathbf{u}^\prime\in\ZZ^n$, $T(\mathbf{u})\cong T(\mathbf{u}^\prime)$ as an $R$-module 
if and only if there exists $y\in\sfM$ such that $u_i=u_i^\prime+\lambda_i(y)$ for all $i=1, \cdots, n$. 
Therefore, we have $\operatorname{Cl}(R)\cong \ZZ^n/\lambda(\ZZ^d)$. 
\end{lemma}

When we can take $\mathbf{u}=\ulcorner \lambda(y)\urcorner$ for some $y\in\RR^d$, a divisorial ideal $T(\mathbf{u})$ is called  \emph{conic}. 
Namely, a divisorial ideal $T(u_1,\cdots,u_n)$ is conic if and only if there exists $y\in\RR^d$ 
such that $u_i-1<\lambda_i(y)\le u_i$ for every $i=1,\cdots,n$. 
Thus, the conic class is also related with hyperplane arrangements and linear programming. 
This nice class of divisorial ideals has been studied in several papers e.g. \cite{Sta,Don,BG1,Bae, Bru}, and they are characterized as follows. 

\begin{proposition}(\cite[Proposition\,3.6]{BG1})
\label{conic}
For a toric singularity $R$, the following is the same set of divisorial ideals. 
\begin{itemize}
\item [(1)] the conic classes, 
\item [(2)] the set of divisorial ideals arising from the decomposition of $R^{1/m}$ as an $R$-module for $m\gg 0$ 
where $R^{1/m}$ is the $R$-module consisted of $m$-th root of elements in $R$. 
\end{itemize}
\end{proposition}

Since $R^{1/m}$ is an MCM $R$-module, a conic divisorial ideal is also an MCM $R$-module. 
Especially, a divisorial ideal which is a torsion element in $\operatorname{Cl}(R)$ is conic \cite[Theorem\,3.2]{BG1}. 
Thus, if $\sigma$ is simplicial then every divisorial ideal is conic. Since $\sigma$ is simplicial if and only if 
every divisorial ideal is an MCM $R$-module \cite[Remark\,4.3]{BG1}, $\operatorname{Cl}(R)$ is larger than conic classes 
when $\sigma$ is not simplicial. 

Further, we remark that if $\mathrm{char}\,k=p$, $R^{1/p}$ is isomorphic to the Frobenius push-forward $F_*R$ 
defined by the Frobenius morphism $F:R\rightarrow R\,\,(r\mapsto r^p)$, 
and the structure of $R^{1/p}$ is important in commutative ring theory in positive characteristic. 
For example, several numerical invariants (e.g. the Hilbert-Kunz multiplicity, the $F$-signature) are computed in \cite[Section~3]{Bru} by paying attention to conic classes.

%%%%%%%%%%%%%%%%%%%%%%%%%%%%%%%%%%%%%%%%%%%%%%%%%%%%%%%%%%%%%%%%%%%%%%%%%%%%%%%%
\subsection{Dimer models and associated quivers}
\label{subsec_dimer}

In this subsection, we introduce a dimer model and define the quiver associated with a dimer model. 
By using this quiver, we will construct the Jacobian algebra, and see such an algebra is an NCCR of a $3$-dimensional Gorenstein toric singularity 
if a dimer model is consistent. 

%\medskip

\begin{figure}[h]
\begin{tikzpicture}

\node (dimer) at (0,0) 
{\scalebox{0.65}{
\begin{tikzpicture}
%vertex
\node (B1) at (1,1){$$}; \node (W1) at (3,3){$$}; 
\draw[thick]  (0,0) rectangle (4,4);

%edge
\draw[line width=0.05cm]  (B1)--(W1); \draw[line width=0.05cm]  (B1)--(0,0);  \draw[line width=0.05cm]  (B1)--(2,0);  \draw[line width=0.05cm]  (B1)--(0,2); 
\draw[line width=0.05cm]  (W1)--(4,4);  \draw[line width=0.05cm]  (W1)--(2,4);  \draw[line width=0.05cm]  (W1)--(4,2); 
%black
\filldraw  [ultra thick, fill=black] (1,1) circle [radius=0.2] ;
%white
\draw  [ultra thick,fill=white] (3,3) circle [radius=0.2] ;
\end{tikzpicture} }};

\node (quiver) at (5,0) 
{\scalebox{0.65}{
\begin{tikzpicture}
%vertex
\draw[lightgray,thick]  (0,0) rectangle (4,4);
\node (B1) at (1,1){$$}; \node (W1) at (3,3){$$}; 
\node (Q0a) at (3.6,0.4){{\Large$0$}}; \node (Q0b) at (0,1){{\Large$0$}}; \node (Q0c) at (3,4){{\Large$0$}}; 
\node (Q1a) at (0.4,3.6){{\Large$1$}}; \node (Q1b) at (1,0){{\Large$1$}}; \node (Q1c) at (4,3){{\Large$1$}}; 

%edge
\draw[lightgray,line width=0.05cm]  (B1)--(W1); \draw[lightgray,line width=0.05cm]  (B1)--(0,0);  
\draw[lightgray,line width=0.05cm]  (B1)--(2,0);  \draw[lightgray,line width=0.05cm]  (B1)--(0,2); 
\draw[lightgray,line width=0.05cm]  (W1)--(4,4);  \draw[lightgray,line width=0.05cm]  (W1)--(2,4);  \draw[lightgray,line width=0.05cm]  (W1)--(4,2); 
%black
\filldraw  [lightgray, ultra thick, fill=lightgray] (1,1) circle [radius=0.2] ;
%white
\draw  [lightgray,ultra thick, fill=white] (3,3) circle [radius=0.2] ;
%arrow
\draw[->, line width=0.08cm] (Q0a)--(Q1a); \draw[->, line width=0.06cm] (Q1b)--(Q0a); \draw[->, line width=0.08cm] (Q1c)--(Q0a);
\draw[->, line width=0.08cm] (Q1a)--(Q0b); \draw[->, line width=0.06cm] (Q1a)--(Q0c); \draw[->, line width=0.08cm] (Q0b)--(Q1b);
\draw[->, line width=0.08cm] (Q0c)--(Q1c);

\end{tikzpicture} }};

\end{tikzpicture} 
\caption{Dimer model and associated quiver}
\label{ex_dimer_A1}
\end{figure}
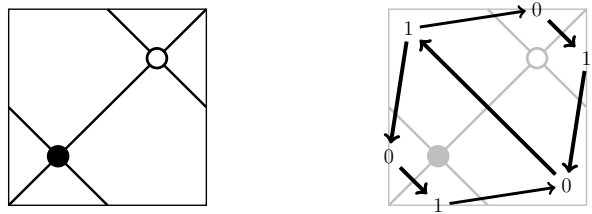

A \emph{dimer model} (or \emph{brane tiling}) is a polygonal cell decomposition of 
the two-torus $\sfT= \RR^2/\ZZ^2$ whose vertices and edges form a finite bipartite graph. 
Thus, each vertex is colored either black or white so that each edge connects a black vertex to a white vertex. 
For example, the left hand side of Figure~\ref{ex_dimer_A1} is a dimer model. 
Let $\Gamma$ be a dimer model. We denote the set of vertices of $\Gamma$ (resp. edges and faces) by $\Gamma_0$ (resp. $\Gamma_1$ and $\Gamma_2$).  

As the dual of a dimer model $\Gamma$, we define the quiver $Q_\Gamma$. 
Sometimes we simply denote it by $Q$ if a situation is clear from the context. 
Namely, we assign a vertex dual to each face in $\Gamma_2$, an arrow dual to each edge in $\Gamma_1$.  
The orientation of arrows is determined so that the white vertex is on the right of the arrow (equivalently, the black vertex is on the left of the arrow). 
For example, the right hand side of Figure~\ref{ex_dimer_A1} is the quiver associated with a dimer model. 
We denote the sets of vertices and arrows of the quiver by $Q_0$ and $Q_1$ respectively. 
In addition, we define the set of oriented faces $Q_2$ as the dual of vertices in a dimer model $\Gamma$. 
The orientation of faces is determined by its boundary. 
That is, by the orientation of each arrow, faces dual to white (resp. black) vertices are oriented clockwise (resp. anti-clockwise). 
%Therefore, we decompose the set of faces as $Q_2=Q^+_2\sqcup Q^-_2$ where $Q^+_2$ (resp. $Q^-_2$) stands for the set of faces oriented clockwise (resp. anti-clockwise). 

\begin{figure}[h] 
\[\scalebox{0.75}{
\begin{tikzpicture}
%vertex
\node (P1) at (0,0){$$}; \node (P2) at (2,0){$$}; 
%edge
\draw[line width=0.035cm]  (P1)--(P2); \draw[line width=0.035cm]  (P1)--(-1.2,1); \draw[line width=0.035cm]  (P1)--(-1.2,-1);
\draw[line width=0.035cm]  (P2)--(3.2,1.2);  \draw[line width=0.035cm]  (P2)--(3.8,0.7); \draw[line width=0.035cm]  (P2)--(3.2,-1); 
%black
\filldraw  [very thick, fill=black] (0,0) circle [radius=0.14] ;
%white
\draw  [very thick, fill=white] (2,0) circle [radius=0.14] ;
%arrow
\draw[->, line width=0.05cm]  (1,-1)--(1,1); \draw[->, line width=0.05cm]  (0.8,1)--(-1.2,0.1); 
%\draw[dotted, line width=0.05cm]  (-1.2,0.4)--(-1.2,-0.4); 
\draw[->, line width=0.05cm]  (-1.2,-0.1)--(0.8,-1); \draw[->, line width=0.05cm]  (1.2,1)--(3.2,0.7); 
\draw[->, line width=0.05cm]  (3.3,0.7)--(3.4, -0.2); 
%\draw[dotted, line width=0.05cm]  (3.6,-0.1)--(3.3,-0.5); 
\draw[->, line width=0.05cm]  (3.4,-0.3)--(1.2,-1); 

\node  at (1.4,0.3) {$a$} ;\node  at (-1.5,0.3) {$p^-_a$} ;\node  at (3.8,0.2) {$p^+_a$} ;
\end{tikzpicture}
}\]
\caption{}
\label{relation}
\end{figure}

Let $\widehat{kQ}$ be the complete path algebra of $Q$. 
For each arrow $a\in Q_1$, there are precisely two oppositely oriented faces which contain the arrow $a$ as a boundary. 
We denote them by $f_a^+, f_a^-\in Q_2$ respectively. Let $p_a^{\pm}$ be the path from $h(a)$ around the boundary of $f^{\pm}_a$ to $t(a)$ (see Figure~\ref{relation}). 
Note that the difference of $p_a^{\pm}$ coincides with the partial derivative $\frac{\partial W_Q}{\partial a}$ of a certain ``potential $W_Q$" 
(for more details, see e.g. \cite[Chapter~2]{Bro}). 
We define the closure of the two-sided ideal $J_Q\coloneqq \langle\overline{p_a^+-p_a^- \mid a\in Q_1}\rangle\subset \widehat{kQ}$ (see \cite[Section~2]{DWZ}). 
Then we define the \emph{complete Jacobian algebra} of a dimer model as $\sfA_Q\coloneqq \widehat{kQ}/J_Q$.

In the rest, we impose the extra condition so-called ``consistency condition". 
Under this assumption, a dimer model will give a non-commutative crepant resolution (see Theorem~\ref{NCCR1}).  

\begin{definition}
\label{def_consistent}
Let $Q$ be a quiver associated with a dimer model $\Gamma$. A map $\sfR:Q_1\rightarrow \RR_{>0}$ which 
satisfies the following conditions is called a \emph{consistent $\sfR$-charge}. 
\begin{itemize}
\item [(1)] $\displaystyle\sum_{a\in\partial f}\sfR(a)=2$ \;for any $f\in Q_2$,  
\item [(2)] $\displaystyle\sum_{h(a)=i}(1-\sfR(a))+\displaystyle\sum_{t(a)=i}(1-\sfR(a))=2$ \;for any $i\in Q_0$. 
\end{itemize}
We say a dimer model $\Gamma$ is \emph{consistent} if it admits a consistent $\sfR$-charge. 
In addition, if a dimer model admits a consistent $\sfR$-charge such that $\sfR(a)<1$ for any $a\in Q_1$, 
we call it an \emph{isoradial} (or \emph{geometrically consistent}) \emph{dimer model}. 
\end{definition} 

By the definition, an isoradial dimer model is consistent. 
In the literature, there are several conditions which are equivalent to the existence of a consistent $\sfR$-charge (see \cite{Boc1, IU1}). 
Note that a consistent $\sfR$-charge is interpreted as $1/\pi$ times of the positive angle arising from the homology classes of some zig-zag paths 
(see e.g. \cite[Chapter~4]{Bro}, \cite[Proof of Theorem~6.6]{Boc1}). 

The complete Jacobian algebra associated with a consistent dimer model has good properties as follows.

\begin{theorem}
\label{NCCR1}
(see e.g. \cite{Bro, IU2, Boc2}) 
Suppose that $Q$ is the quiver associated with a consistent dimer model and $\sfA_Q$ is the complete Jacobian algebra. 
Let $R\coloneqq \rmZ(\sfA_Q)$ be the center of $\sfA_Q$. 
Then $R$ is a $3$-dimensional complete local Gorenstein toric singularity, and 
there exists a generator $M\in\CM R$ such that $\sfA_Q\cong\End_R(M)$ as an $R$-algebra and $M$ is a direct sum of rank one MCM modules. 
Moreover, $\End_R(M)$ is an NCCR of $R$. 
\end{theorem} 

\begin{remark}
\label{toric_diagram}
We note that it is known that $R$ is Gorenstein if and only if there is $x\in\sigma^\vee\cap\ZZ^3$ such that $\lambda_i(x)=1$ for all $i=1, \cdots, n$ 
where $v_1, \cdots, v_n\in\ZZ^3$ are minimal generators of $\sigma$ (see e.g. \cite[Theorem~6.33]{BG2}). 
Thus, for a given $3$-dimensional Gorensitein toric singularity, 
we can take a hyperplane $z=1$ so that generators $v_1,\cdots,v_n$ lie on this hyperplane (i.e. the third coordinate of $v_i$ is $1$). 
Thus, we obtain the lattice polygon $\Delta\subset\RR^2$ on this plane. 
Conversely, for a lattice polygon $\Delta$ in $\RR^2$, we define the cone $\sigma_\Delta\subset\RR^3$ whose section on the hyperplane $z=1$ is $\Delta$. 
Then the associated toric singularity is Gorenstein in dimension three. 
\end{remark}

In this way, we obtain a $3$-dimensional complete local Gorenstein toric singularity $R=k[[\sigma^\vee\cap\ZZ^3]]$ and its NCCR from a consistent dimer model. 
Also, for every lattice polygon $\Delta$ in $\RR^2$ (equivalently, for every $3$-dimensional Gorenstein toric singularity $R=k[[\sigma_\Delta^\vee\cap\ZZ^3]]$), 
there exists a consistent dimer model giving $R$ as the center of the complete Jacobian algebra (see \cite{Gul, IU2}). 
Thus, every $3$-dimensional Gorenstein toric singularity admits an NCCR which is constructed from a consistent dimer model. 
%Note that consistent dimer models arising from the method in \cite{Gul} are isoradial. 
%On the other hand, ones arising from the method in \cite{IU2} are not necessarily isoradial. 
In general, a consistent dimer model which gives an NCCR of $R$ is not unique. 
As we showed in Theorem~\ref{NCCR1}, a generator $M\in\CM R$ which satisfies $\sfA_Q\cong\End_R(M)$ is a finite direct sum of rank one MCM modules. 
Conversely, Bocklandt showed every $R$-module $M$ giving an NCCR of $R$ is always coming from a consistent dimer model if $M$ is a finite direct sum of rank one MCM $R$-modules \cite{Boc3}.

%%%%%%%%%%%%%%%%%%%%%%%%%%%%%%%%%%%%%%%%%%%%%%%%%%%%%%%%%%%%
%%%%%%%%%%%%%%%%%%%%%%%%%%%%%%%%%%%%%%%%%%%%%%%%%%%%%%%%%%%%
%%%%%%%%%%%%%%%%%%%%%%%%%%%%%%%%%%%%%%%%%%%%%%%%%%%%%%%%%%%%
\section{Conic divisorial ideals arising from dimer models} 
\label{sec_conic_dimer}

\subsection{Examples}
We start this section with observing some examples. 
In the following examples, we can see that all conic divisorial ideals arise from an isoradial dimer model. 

Firstly, we consider the case of the $A_1$-singularity. 

%%%%%%%%%%%%%%%%%%%%%%%%%%%%%%%%%%%%%%%%%%%%%%%%%%%%%%%%%%%%
\begin{example}
\label{ex_conifold}
Let $R=k[[x, y, z, w]]/(xy-zw)$ be a $3$-dimensional $A_1$-singularity (or conifold singularity). 
$R$ is a toric singularity, and further it is a simple singularity. 
Therefore, $R$ is of finite CM representation type, that is, it has only finitely many non-isomorphic indecomposable MCM modules. 
It is known that finitely many MCM modules are $R$, $T(1,0,0,0)$ and $T(0,1,0,0)\cong T(1,0,0,0)^*$. 
Also, modules which give an NCCR of $R$ are only $M_1\coloneqq R\oplus T(1,0,0,0)$ and $M_2\coloneqq R\oplus T(0,1,0,0)$ (see e.g. \cite{VdB2}), 
and it is easy to see $\End_R(M_1)\cong\End_R(M_2)$ as an $R$-algebra. 
Note that the consistent dimer model $\Gamma$ giving this endomorphism ring as the complete Jacobian algebra is just Figure~\ref{ex_dimer_A1}. 
We can easily see that this dimer model is isoradial. 
Let $\sfA$ be the complete Jacobian algebra of $\Gamma$. Then we have 
\[
\add_R(\sfA)=\add_R(\End_R(M_1))=\add_R(R\oplus T(1,0,0,0)\oplus T(0,1,0,0)). 
\]
On the other hand, we can check all MCM modules are conic by a direct computation (see also \cite[Proof of Theorem~6.1]{TY}). 
Thus, for $m\gg 0$, we have 
\[\add_R(R^{1/m})=\add_R(R\oplus T(1,0,0,0)\oplus T(0,1,0,0))\]. 
\end{example}

%%%%%%%%%%%%%%%%%%%%%%%%%%%%%%%%%%%%%%%%%%%%%%%%%%%%%%%%%%%% 
Next, we show simplicial cases. 

\begin{example}
\label{ex_simplicial}
Let $\Delta$ be a triangle polygon in $\RR^2$, then the cone $\sigma_\Delta\subset\RR^3$ is simplicial. 
The associated toric singularity $R$ is a $3$-dimensional Gorenstein quotient singularity by a finite abelian group (see e.g. \cite[Example~1.3.20]{CLS}). 
Thus, we set $R=S^G$ where $S=k[[x, y, z]]$ and $G\subset\SL(3,k)$ is a finite abelian group with $(\mathrm{char}\,k,|G|)=1$, and we may assume $G$ is small. 
Let $V_0=k, V_1, \cdots, V_{|G|-1}$ be the full set of non-isomorphic irreducible representations of $G$, and set $M_i\coloneqq(S\otimes_kV_i)^G$. 
Then these give all rank one MCM $R$-modules. 
Furthermore, the $R$-module $S$ is decomposed as $S\cong R\oplus M_1\oplus\cdots\oplus M_{|G|-1}$, 
and this gives NCCR of $R$ (see e.g. \cite[Example~2.3]{IN}). 
Also, a consistent dimer model associated with $R$ is homotopy equivalent to a regular hexagonal dimer model 
(i.e. each face of a dimer model is a regular hexagon).
We denote such a dimer model by $\Gamma$, and the associated quiver by $Q$. 
By defining an $\sfR$-charge as $\sfR(a)=\frac{2}{3}$ for any $a\in Q_1$, we see that $\Gamma$ is isoradial. 
Also, we note that $Q$ is the McKay quiver of $G$, and we have $\sfA_Q\cong\End_R(S)\cong S*G$ (see \cite{UY}, \cite[1.7, 1.8]{IN}).  
Thus, we have 
\[
\add_R(\sfA_Q)=\add_R(S)=\add_R(R\oplus M_1\oplus\cdots\oplus M_{|G|-1}).
\] 
Since $\{[R], [M_1],\cdots, [M_{|G|-1}]\}$ are torsion elements in $\operatorname{Cl}(R)$, they are all conic. 
\end{example}

%%%%%%%%%%%%%%%%%%%%%%%%%%%%%%%%%%%%%%%%%%%%%%%%%%%%%%%%%%%%
\subsection{Toric singularities associated with reflexive polygons}
\label{subsec_reflexive}

In Example~\ref{ex_conifold} and \ref{ex_simplicial}, we saw that conic divisorial ideals arise from an isoradial dimer model. 
In this subsection, we observe this phenomenon for the case of toric singularities associated with reflexive polygons, 
and prove Theorem~\ref{main} and Corollary~\ref{main_cor}.

\begin{figure}[h]
\begin{center}
\scalebox{0.75}{
\begin{tikzpicture}

\node (3a) at (0,0){
%%% 3a %%%
\scalebox{0.65}{
\begin{tikzpicture}
\draw [step=1,thin, gray] (-2,-2) grid (2,2);
\filldraw [ultra thick, fill=white] (1,0)--(0,1)--(-1,-1)--cycle ;
\filldraw  [ultra thick, fill=black] (1,0) circle [radius=0.1] ;
\filldraw  [ultra thick, fill=black] (0,1) circle [radius=0.1] ;
\filldraw  [ultra thick, fill=black] (-1,-1) circle [radius=0.1] ;

\node  at (0,0) {\LARGE{3a}} ;
\end{tikzpicture} }} ;

%%% 4a %%%
\node (4a) at (4,0){
\scalebox{0.65}{
\begin{tikzpicture}
\draw [step=1,thin, gray] (-2,-2) grid (2,2);
\filldraw [ultra thick, fill=white] (1,0)--(0,1)--(-1,0)--(0,-1)--cycle ;
\filldraw  [ultra thick, fill=black] (1,0) circle [radius=0.1] ;
\filldraw  [ultra thick, fill=black] (0,1) circle [radius=0.1] ;
\filldraw  [ultra thick, fill=black] (-1,0) circle [radius=0.1] ;
\filldraw  [ultra thick, fill=black] (0,-1) circle [radius=0.1] ;

\node  at (0,0) {\LARGE{4a}} ;
\end{tikzpicture} }};

%%% 4b %%%
\node (4b) at (8,0){
\scalebox{0.65}{
\begin{tikzpicture}
\draw [step=1,thin, gray] (-2,-2) grid (2,2);
\filldraw [ultra thick, fill=white] (0,1)--(-1,0)--(0,-1)--(1,-1)--cycle ;
\filldraw  [ultra thick, fill=black] (0,1) circle [radius=0.1] ;
\filldraw  [ultra thick, fill=black] (-1,0) circle [radius=0.1] ;
\filldraw  [ultra thick, fill=black] (0,-1) circle [radius=0.1] ;
\filldraw  [ultra thick, fill=black] (1,-1) circle [radius=0.1] ;

\node  at (0,0) {\LARGE{4b}} ;
\end{tikzpicture} }}; 

%%% 4c %%%
\node (4c) at (12,0){
\scalebox{0.65}{
\begin{tikzpicture}
\draw [step=1,thin, gray] (-2,-2) grid (2,2);
\filldraw [ultra thick, fill=white] (0,1)--(-1,-1)--(1,-1)--cycle ;
\filldraw  [ultra thick, fill=black] (0,1) circle [radius=0.1] ;
\filldraw  [ultra thick, fill=black] (-1,-1) circle [radius=0.1] ;
\filldraw  [ultra thick, fill=black] (1,-1) circle [radius=0.1] ;

\node  at (0,0) {\LARGE{4c}} ;
\end{tikzpicture} }};

%%% 5a %%%
\node (5a) at (0,-3){
\scalebox{0.65}{
\begin{tikzpicture}
\draw [step=1, thin, gray] (-2,-2) grid (2,2);
\filldraw [ultra thick, fill=white] (1,0)--(0,1)--(-1,1)--(-1,0)--(0,-1)--cycle ;
\filldraw  [ultra thick, fill=black] (1,0) circle [radius=0.1] ;
\filldraw  [ultra thick, fill=black] (0,1) circle [radius=0.1] ;
\filldraw  [ultra thick, fill=black] (-1,1) circle [radius=0.1] ;
\filldraw  [ultra thick, fill=black] (-1,0) circle [radius=0.1] ;
\filldraw  [ultra thick, fill=black] (0,-1) circle [radius=0.1] ;

\node  at (0,0) {\LARGE{5a}} ;
\end{tikzpicture} }}; 

%%% 5b %%%
\node (5b) at (4,-3){
\scalebox{0.65}{
\begin{tikzpicture}
\draw [step=1,thin, gray] (-2,-2) grid (2,2);
\filldraw [ultra thick, fill=white] (0,1)--(-1,0)--(-1,-1)--(1,-1)--cycle ;
\filldraw  [ultra thick, fill=black] (0,1) circle [radius=0.1] ;
\filldraw  [ultra thick, fill=black] (-1,0) circle [radius=0.1] ;
\filldraw  [ultra thick, fill=black] (-1,-1) circle [radius=0.1] ;
\filldraw  [ultra thick, fill=black] (0,-1) circle [radius=0.1] ;
\filldraw  [ultra thick, fill=black] (1,-1) circle [radius=0.1] ;

\node  at (0,0) {\LARGE{5b}} ;
\end{tikzpicture} }};

%%% 6a %%%
\node (6a) at (8,-3){
\scalebox{0.65}{
\begin{tikzpicture}
\draw [step=1,thin, gray] (-2,-2) grid (2,2);
\filldraw [ultra thick, fill=white] (1,0)--(0,1)--(-1,1)--(-1,0)--(0,-1)--(1,-1)--cycle ;
\filldraw  [ultra thick, fill=black] (1,0) circle [radius=0.1] ;
\filldraw  [ultra thick, fill=black] (0,1) circle [radius=0.1] ;
\filldraw  [ultra thick, fill=black] (-1,1) circle [radius=0.1] ;
\filldraw  [ultra thick, fill=black] (-1,0) circle [radius=0.1] ;
\filldraw  [ultra thick, fill=black] (0,-1) circle [radius=0.1] ;
\filldraw  [ultra thick, fill=black] (1,-1) circle [radius=0.1] ;

\node  at (0,0) {\LARGE{6a}} ;
\end{tikzpicture} }};

%%% 6b %%%
\node (6b) at (12,-3){
\scalebox{0.65}{
\begin{tikzpicture}
\draw [step=1,thin, gray] (-2,-2) grid (2,2);
\filldraw [ultra thick, fill=white] (1,0)--(0,1)--(-1,1)--(-1,-1)--(0,-1)--cycle ;
\filldraw  [ultra thick, fill=black] (1,0) circle [radius=0.1] ;
\filldraw  [ultra thick, fill=black] (0,1) circle [radius=0.1] ;
\filldraw  [ultra thick, fill=black] (-1,1) circle [radius=0.1] ;
\filldraw  [ultra thick, fill=black] (-1,0) circle [radius=0.1] ;
\filldraw  [ultra thick, fill=black] (-1,-1) circle [radius=0.1] ;
\filldraw  [ultra thick, fill=black] (0,-1) circle [radius=0.1] ;

\node  at (0,0) {\LARGE{6b}} ;
\end{tikzpicture} }};

%%% 6c %%%
\node (6c) at (0,-6){
\scalebox{0.65}{
\begin{tikzpicture}
\draw [step=1,thin, gray] (-2,-2) grid (2,2);
\filldraw [ultra thick, fill=white] (1,0)--(0,1)--(-2,-1)--(0,-1)--cycle ;
\filldraw  [ultra thick, fill=black] (1,0) circle [radius=0.1] ;
\filldraw  [ultra thick, fill=black] (0,1) circle [radius=0.1] ;
\filldraw  [ultra thick, fill=black] (-1,0) circle [radius=0.1] ;
\filldraw  [ultra thick, fill=black] (-2,-1) circle [radius=0.1] ;
\filldraw  [ultra thick, fill=black] (-1,-1) circle [radius=0.1] ;
\filldraw  [ultra thick, fill=black] (0,-1) circle [radius=0.1] ;

\node  at (0,0) {\LARGE{6c}} ;
\end{tikzpicture} }};

%%% 6d %%%
\node (6d) at (4,-6){
\scalebox{0.65}{
\begin{tikzpicture}
\draw [step=1,thin, gray] (-2,-2) grid (2,2);
\filldraw [ultra thick, fill=white] (0,1)--(-2,-1)--(1,-1)--cycle ;
\filldraw  [ultra thick, fill=black] (0,1) circle [radius=0.1] ;
\filldraw  [ultra thick, fill=black] (-1,0) circle [radius=0.1] ;
\filldraw  [ultra thick, fill=black] (-2,-1) circle [radius=0.1] ;
\filldraw  [ultra thick, fill=black] (-1,-1) circle [radius=0.1] ;
\filldraw  [ultra thick, fill=black] (0,-1) circle [radius=0.1] ;
\filldraw  [ultra thick, fill=black] (1,-1) circle [radius=0.1] ;

\node  at (0,0) {\LARGE{6d}} ;
\end{tikzpicture} }};

%%% 7a %%%
\node (7a) at (8,-6){
\scalebox{0.65}{
\begin{tikzpicture}
\draw [step=1,thin, gray] (-2,-2) grid (2,2);
\filldraw [ultra thick, fill=white] (1,0)--(0,1)--(-1,1)--(-1,-1)--(1,-1)--cycle ;
\filldraw  [ultra thick, fill=black] (1,0) circle [radius=0.1] ;
\filldraw  [ultra thick, fill=black] (0,1) circle [radius=0.1] ;
\filldraw  [ultra thick, fill=black] (-1,1) circle [radius=0.1] ;
\filldraw  [ultra thick, fill=black] (-1,0) circle [radius=0.1] ;
\filldraw  [ultra thick, fill=black] (-1,-1) circle [radius=0.1] ;
\filldraw  [ultra thick, fill=black] (0,-1) circle [radius=0.1] ;
\filldraw  [ultra thick, fill=black] (1,-1) circle [radius=0.1] ;

\node  at (0,0) {\LARGE{7a}} ;
\end{tikzpicture} }};

%%% 7b %%%
\node (7b) at (12,-6){
\scalebox{0.65}{
\begin{tikzpicture}
\draw [step=1,thin, gray] (-2,-2) grid (2,2);
\filldraw [ultra thick, fill=white] (1,0)--(0,1)--(-2,-1)--(1,-1)--cycle ;
\filldraw  [ultra thick, fill=black] (1,0) circle [radius=0.1] ;
\filldraw  [ultra thick, fill=black] (0,1) circle [radius=0.1] ;
\filldraw  [ultra thick, fill=black] (-1,0) circle [radius=0.1] ;
\filldraw  [ultra thick, fill=black] (-2,-1) circle [radius=0.1] ;
\filldraw  [ultra thick, fill=black] (-1,-1) circle [radius=0.1] ;
\filldraw  [ultra thick, fill=black] (0,-1) circle [radius=0.1] ;
\filldraw  [ultra thick, fill=black] (1,-1) circle [radius=0.1] ;

\node  at (0,0) {\LARGE{7b}} ;
\end{tikzpicture} }};

%%% 8a %%%
\node (8a) at (0,-9){
\scalebox{0.65}{
\begin{tikzpicture}
\draw [step=1,thin, gray] (-2,-2) grid (2,2);
\filldraw [ultra thick, fill=white] (1,1)--(-1,1)--(-1,-1)--(1,-1)--cycle ;
\filldraw  [ultra thick, fill=black] (1,0) circle [radius=0.1] ;
\filldraw  [ultra thick, fill=black] (1,1) circle [radius=0.1] ;
\filldraw  [ultra thick, fill=black] (0,1) circle [radius=0.1] ;
\filldraw  [ultra thick, fill=black] (-1,1) circle [radius=0.1] ;
\filldraw  [ultra thick, fill=black] (-1,0) circle [radius=0.1] ;
\filldraw  [ultra thick, fill=black] (-1,-1) circle [radius=0.1] ;
\filldraw  [ultra thick, fill=black] (0,-1) circle [radius=0.1] ;
\filldraw  [ultra thick, fill=black] (1,-1) circle [radius=0.1] ;

\node  at (0,0) {\LARGE{8a}} ;
\end{tikzpicture} }};

%%% 8b %%%
\node (8b) at (4,-9){
\scalebox{0.65}{
\begin{tikzpicture}
\draw [step=1,thin, gray] (-2,-2) grid (2,2);
\filldraw [ultra thick, fill=white] (0,1)--(-1,1)--(-1,-1)--(2,-1)--cycle ;
\filldraw  [ultra thick, fill=black] (1,0) circle [radius=0.1] ;
\filldraw  [ultra thick, fill=black] (0,1) circle [radius=0.1] ;
\filldraw  [ultra thick, fill=black] (-1,1) circle [radius=0.1] ;
\filldraw  [ultra thick, fill=black] (-1,0) circle [radius=0.1] ;
\filldraw  [ultra thick, fill=black] (-1,-1) circle [radius=0.1] ;
\filldraw  [ultra thick, fill=black] (0,-1) circle [radius=0.1] ;
\filldraw  [ultra thick, fill=black] (1,-1) circle [radius=0.1] ;
\filldraw  [ultra thick, fill=black] (2,-1) circle [radius=0.1] ;

\node  at (0,0) {\LARGE{8b}} ;
\end{tikzpicture} }};

%%% 8c %%%
\node (8c) at (8,-9){
\scalebox{0.65}{
\begin{tikzpicture}
\draw [step=1,thin, gray] (-2,-2) grid (2,2);
\filldraw [ultra thick, fill=white] (0,1)--(-2,-1)--(2,-1)--cycle ;
\filldraw  [ultra thick, fill=black] (1,0) circle [radius=0.1] ;
\filldraw  [ultra thick, fill=black] (0,1) circle [radius=0.1] ;
\filldraw  [ultra thick, fill=black] (-1,0) circle [radius=0.1] ;
\filldraw  [ultra thick, fill=black] (-2,-1) circle [radius=0.1] ;
\filldraw  [ultra thick, fill=black] (-1,-1) circle [radius=0.1] ;
\filldraw  [ultra thick, fill=black] (0,-1) circle [radius=0.1] ;
\filldraw  [ultra thick, fill=black] (1,-1) circle [radius=0.1] ;
\filldraw  [ultra thick, fill=black] (2,-1) circle [radius=0.1] ;

\node  at (0,0) {\LARGE{8c}} ;
\end{tikzpicture} }};

%%% 9a %%%
\node (9a) at (12,-9){
\scalebox{0.65}{
\begin{tikzpicture}
\draw [step=1,thin, gray] (-2,-2) grid (2,2);
\filldraw [ultra thick, fill=white] (-1,2)--(-1,-1)--(2,-1)--cycle ;
\filldraw  [ultra thick, fill=black] (1,0) circle [radius=0.1] ;
\filldraw  [ultra thick, fill=black] (0,1) circle [radius=0.1] ;
\filldraw  [ultra thick, fill=black] (-1,2) circle [radius=0.1] ;
\filldraw  [ultra thick, fill=black] (-1,1) circle [radius=0.1] ;
\filldraw  [ultra thick, fill=black] (-1,0) circle [radius=0.1] ;
\filldraw  [ultra thick, fill=black] (-1,-1) circle [radius=0.1] ;
\filldraw  [ultra thick, fill=black] (0,-1) circle [radius=0.1] ;
\filldraw  [ultra thick, fill=black] (1,-1) circle [radius=0.1] ;
\filldraw  [ultra thick, fill=black] (2,-1) circle [radius=0.1] ;

\node  at (0,0) {\LARGE{9a}} ;
\end{tikzpicture} }};

\end{tikzpicture} }
\end{center}
\caption{The classification of reflexive polygons}
\label{class_ref}
\end{figure}
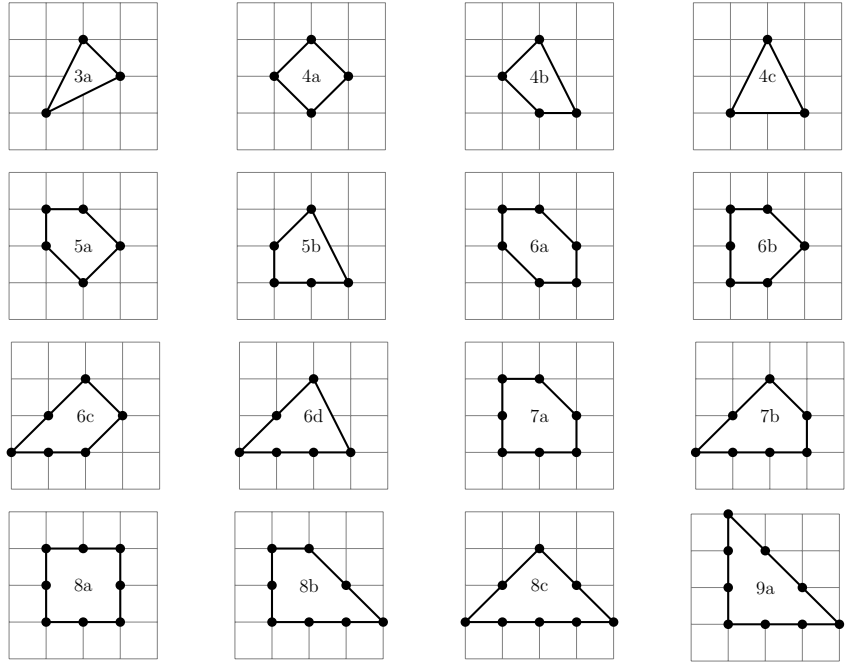

In what follows, we consider a $3$-dimensional Gorenstein toric singularity whose associated lattice polygon $\Delta$ is a reflexive polygon. 
We recall that $\Delta$ is called a reflexive polygon (or Fano polygon) if the origin is the unique interior point of $\Delta$. 
Reflexive polygons are classified in $16$ types (see Figure~\ref{class_ref}) up to integral unimodular transformations 
(see e.g. \cite[Theorem~8.3.7]{CLS}, \cite[Appendix]{Boc2}).  

%%%%%%%%%%%%%%%%%%%%%%%%%%%%%%%%%%%%%%%%%%%%%%%%%%%%%%%

For these cases, the associated consistent dimer models are well-studied in e.g. \cite{Boc2, Boc4, HS, Nak}, 
and such dimer models are completely classified. 
For a given dimer model, we can see whether it is isoradial or not by checking Definition~\ref{def_consistent} or an equivalent condition as in \cite[Theorem~5.1]{KS}.
Also, for these cases, we already know the precise description of $M_i$'s appearing in the decomposition (\ref{decomp_NCCR}) by the results of \cite[Section~5]{Nak}. 
(Note that $\sfA_Q\cong\bigoplus_{i,j\in Q_0} T_{ij}$ as an $R$-module by the notation in \cite{Nak}.) 
Therefore, we prove Theorem~\ref{main} and Corollary~\ref{main_cor} by a case-by-case check. 
Note that a divisorial ideal $I$ is conic if and only if the $R$-dual $I^*$ is conic (see \cite[Remark~3.4 (b)]{BG1}). 
Thus, to check a given divisorial ideal is conic or not, we may only consider one of them. 
Also, the next lemma is useful to check the conicness. 

\begin{lemma} 
\label{conic_check}
(see \cite[Proposition~1.4]{Bru}) 
A divisorial ideal $T(\mathbf{u})=T(u_1,\cdots,u_n)$ is conic if and only if there is $y\in\RR^d$ such that $\lambda(y)-\mathbf{u}$ 
is in the semi-open cube $(-1,0]^n$. 
\end{lemma}

\medskip

Now, we consider Theorem~\ref{main} and Corollary~\ref{main_cor}.  
Remark that type 3a, 4c, 6d, 8c and 9a are contained in Example~\ref{ex_simplicial}, because the associated cones are simplicial. 
Thus, we investigate other cases. 

%%%%%%%%%%%%%%%%%%%%%%%%%%%%%%%%%%%%%%%%%%%%%%%%%%%%%%%
\subsubsection{{\bf \textup{Type 4a}}}
\label{type4a}

We consider the reflexive polygon of type 4a. 
Thus, let $R$ be the $3$-dimensional complete local Gorenstein toric singularity defined by the cone $\sigma$: 
\[
\sigma=\mathrm{Cone}\{v_1=(1,0,1), v_2=(0,1,1), v_3=(-1,0,1), v_4=(0,-1,1) \}. 
\]
Every divisorial ideal (i.e. rank one reflexive $R$-module) takes the form $T(a, b, c, d)$, and it corresponds to  
the element $a[I_1]+\cdots+d[I_4]$ in $\Cl(R)$ (see subsection~\ref{toric_pre}). 
As elements in $\Cl(R)$, we have $[I_1]=[I_3]$, $[I_2]=[I_4]$, $2[I_1]+2[I_2]=0$ (see Lemma~\ref{div_eq}). 
Thus, we have $\Cl(R)\cong\ZZ\times\ZZ/2\ZZ$, and each divisorial ideal is represented by $T(a,b,0,0)$ where $a\in\ZZ, b\in\ZZ/2\ZZ$. 
In this case, there are two consistent dimer models associated with $R$.   
By the results in \cite[subsection~5.2]{Nak}, we have rank one MCM $R$-modules as in Figure~\ref{4a} from these consistent dimer models. 
In this figure, each circle stands for $(a,b)\in\Cl(R)$ corresponding to an MCM module $T(a,b,0,0)$. 
Especially, a double circle stands for the origin $(0,0)$.
Also, by the direct computation below, we can see black circles are conic.  
 
\begin{center}
\begin{tabular}{ll}
$T(0,1,0,0)\cong T(\lambda(\frac{1}{4},\frac{1}{2},\frac{-1}{4}))$,& $T(1,0,0,0)\cong T(\lambda(\frac{1}{2},\frac{1}{4},\frac{-1}{4}))$,\\  
$T(1,1,0,0)\cong T(\lambda(\frac{1}{2},\frac{1}{2},0))$. 
\end{tabular} 
\end{center}
On the other hand, we see that white circles are not conic. 
For example, if $T(3,1,0,0)$ is conic, then there is $(\alpha,\beta,\gamma)\in\RR^3$ such that $\lambda(\alpha,\beta,\gamma)-(3,1,0,0)\in(-1,0]^4$ 
by Lemma~\ref{conic_check}. Thus, we can write 
\[
(\alpha+\gamma-3,\beta+\gamma-1,-\alpha+\gamma,-\beta+\gamma)=(\epsilon_1,\epsilon_2,\epsilon_3,\epsilon_4) 
\]
for some $\epsilon_i\in(-1,0]$. Then we easily have $\epsilon_1+\epsilon_3+2=\epsilon_2+\epsilon_4$.  
Here, the left hand side is greater than $0$, while the right hand side is less than or equal to $0$. 
Hence we have the contradiction. 

\begin{figure}[h]
{\scalebox{0.6}{
\begin{tikzpicture}
\draw [step=1,thin, gray] (-2,-1) grid (4,2);

\filldraw  [thick, fill=black] (-1,0) circle [radius=0.17] ;
\draw (0,0) circle [radius=0.3]; 
\filldraw (0,0) circle [radius=0.15]; 
\filldraw  [thick, fill=black] (0,1) circle [radius=0.17] ;
\filldraw  [thick, fill=black] (1,0) circle [radius=0.17] ;
\filldraw  [thick, fill=black] (1,1) circle [radius=0.17] ;
\filldraw  [thick, fill=black] (2,1) circle [radius=0.17] ;

\filldraw  [thick, fill=white] (-1,1) circle [radius=0.17] ;
\filldraw  [thick, fill=white] (3,1) circle [radius=0.17] ;

%\node  at (1,-1.8) {{\Large\text{MCMs arising from consistent dimer models}}} ;
\end{tikzpicture}
}}
\caption{MCMs arising from consistent dimer models for Type 4a} 
\label{4a}
\end{figure}
 
As we mentioned, there are two consistent dimer models associated with $R$, 
and one of them is isoradial and the other is not isoradial. 
The following figures are an isoradial one and the associated quiver. (The fraction on each arrow is an $\sfR$-charge which indicates this dimer model is isoradial.) 

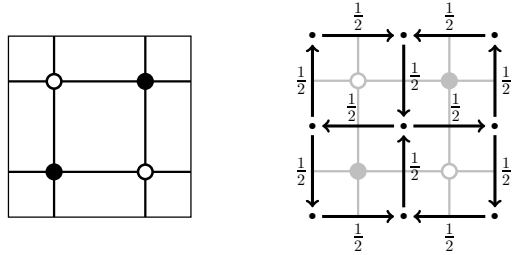
\begin{center}
\begin{tikzpicture} 
\node (DM) at (0,0) 
{\scalebox{0.6}{
\begin{tikzpicture}
%vertex
\node (P1) at (1,1){$$}; \node (P2) at (3,1){$$}; \node (P3) at
(3,3){$$}; \node (P4) at (1,3){$$};

\draw[thick]  (0,0) rectangle (4,4);

%edge
\draw[line width=0.05cm]  (P1)--(P2)--(P3)--(P4)--(P1);\draw[line width=0.05cm] (0,1)--(P1)--(1,0); \draw[line width=0.05cm]  (4,1)--(P2)--(3,0);
\draw[line width=0.05cm]  (0,3)--(P4)--(1,4);\draw[line width=0.05cm]  (3,4)--(P3)--(4,3);
%black
\filldraw  [ultra thick, fill=black] (1,1) circle [radius=0.16] ;\filldraw  [ultra thick, fill=black] (3,3) circle [radius=0.16] ;
%white
\draw  [ultra thick,fill=white] (3,1) circle [radius=0.16] ;\draw  [ultra thick, fill=white] (1,3) circle [radius=0.16] ;
\end{tikzpicture} 
} } ;

\node (Qui) at (4,0) 
{\scalebox{0.6}{
\begin{tikzpicture}
%vertex
\node (Q1) at (2,2){$\bullet$};\node (Q2a) at (0,2){$\bullet$}; \node(Q2b) at (4,2){$\bullet$};\node (Q3a) at (0,0){$\bullet$};
\node(Q3c) at (4,4){$\bullet$};\node(Q3b) at (4,0){$\bullet$};\node(Q3d) at (0,4){$\bullet$};\node (Q4a) at (2,0){$\bullet$};
\node (Q4b) at (2,4){$\bullet$};

%\node (Q1) at (2,2){$0$};\node (Q2a) at (0,2){$1$}; \node(Q2b) at (4,2){$1$};\node (Q3a) at (0,0){$2$};
%\node(Q3c) at (4,4){$2$};\node(Q3b) at (4,0){$2$};\node(Q3d) at (0,4){$2$};\node (Q4a) at (2,0){$3$};
%\node (Q4b) at (2,4){$3$};

%edge
\draw[lightgray, line width=0.05cm]  (P1)--(P2)--(P3)--(P4)--(P1);\draw[lightgray, line width=0.05cm] (0,1)--(P1)--(1,0); 
\draw[lightgray, line width=0.05cm]  (4,1)--(P2)--(3,0);\draw[lightgray, line width=0.05cm]  (0,3)--(P4)--(1,4);\draw[lightgray, line width=0.05cm]  (3,4)--(P3)--(4,3);
%black
\filldraw  [ultra thick, draw=lightgray, fill=lightgray] (1,1) circle [radius=0.16] ;\filldraw  [ultra thick, draw=lightgray, fill=lightgray] (3,3) circle [radius=0.16] ;
%white
\draw  [ultra thick, draw=lightgray,fill=white] (3,1) circle [radius=0.16] ;\draw  [ultra thick, draw=lightgray,fill=white] (1,3) circle [radius=0.16] ;

\draw[->, line width=0.067cm] (Q1)--node[midway,xshift=-4pt,yshift=12pt] {{\LARGE$\frac{1}{2}$}}(Q2a);
\draw[->, line width=0.067cm] (Q2a)--node[midway,xshift=-7pt] {{\LARGE$\frac{1}{2}$}}(Q3a);
\draw[->, line width=0.067cm] (Q3a)--node[below,midway,yshift=0pt] {{\LARGE$\frac{1}{2}$}}(Q4a);
\draw[->, line width=0.067cm] (Q4a)--node[midway,xshift=7pt,yshift=3pt] {{\LARGE$\frac{1}{2}$}}(Q1);
\draw[->, line width=0.067cm] (Q2a)--node[midway,xshift=-7pt] {{\LARGE$\frac{1}{2}$}}(Q3d);
\draw[->, line width=0.067cm] (Q3d)--node[midway,yshift=12pt] {{\LARGE$\frac{1}{2}$}}(Q4b);
\draw[->, line width=0.067cm] (Q4b)--node[midway,xshift=7pt,yshift=3pt] {{\LARGE$\frac{1}{2}$}}(Q1);
\draw[->, line width=0.067cm] (Q1)--node[midway,xshift=4pt,yshift=12pt] {{\LARGE$\frac{1}{2}$}}(Q2b);
\draw[->, line width=0.067cm] (Q2b)--node[midway,xshift=7pt] {{\LARGE$\frac{1}{2}$}}(Q3b);
\draw[->, line width=0.067cm] (Q3b)--node[below,midway,yshift=0pt] {{\LARGE$\frac{1}{2}$}}(Q4a);
\draw[->, line width=0.067cm] (Q2b)--node[midway,xshift=7pt] {{\LARGE$\frac{1}{2}$}}(Q3c);
\draw[->, line width=0.067cm] (Q3c)--node[midway,yshift=12pt] {{\LARGE$\frac{1}{2}$}}(Q4b);
\end{tikzpicture} 
} } ;
\end{tikzpicture}
\end{center}
Again, by the results in \cite[subsection~5.2]{Nak}, we can see rank one MCM modules arising from this isoradial dimer model are actually conic, 
and non-conic divisorial ideals arise from a consistent dimer model which is not isoradial.

%%%%%%%%%%%%%%%%%%%%%%%%%%%%%%%%%%%%%%%%%%%%%%%%%%%%%%%
\subsubsection{{\bf \textup{Type 4b}}}
\label{type4b}

\begin{figure}[h]
{\scalebox{0.6}{
\begin{tikzpicture}
\draw [step=1,thin, gray] (-4,-1) grid (4,1);

\draw (0,0) circle [radius=0.3]; 
\filldraw (0,0) circle [radius=0.15]; 
\filldraw  [thick, fill=black] (1,0) circle [radius=0.17] ;
\filldraw  [thick, fill=black] (2,0) circle [radius=0.17] ;
\filldraw  [thick, fill=black] (3,0) circle [radius=0.17] ;
\filldraw  [thick, fill=black] (-1,0) circle [radius=0.17] ;
\filldraw  [thick, fill=black] (-2,0) circle [radius=0.17] ;
\filldraw  [thick, fill=black] (-3,0) circle [radius=0.17] ;

%\node  at (1,-1.8) {{\Large\text{MCMs arising from consistent dimer models}}} ;
\end{tikzpicture}
}}
\caption{MCMs arising from consistent dimer models for Type 4b} 
\label{4b}
\end{figure}
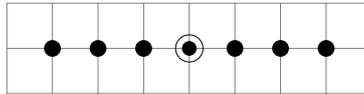

We consider the reflexive polygon of type 4b. 
Thus, let $R$ be the $3$-dimensional complete local Gorenstein toric singularity defined by the cone $\sigma$: 
\[
\sigma=\mathrm{Cone}\{v_1=(0,1,1), v_2=(-1,0,1), v_3=(0,-1,1), v_4=(1,-1,1) \}. 
\]
Then we consider a divisorial ideal $T(a,b,c,d)$. 
As elements in $\Cl(R)$, we have $[I_3]=3[I_1]$, $2[I_1]+[I_2]=0$, $[I_2]=[I_4]$.  
Therefore, we have $\Cl(R)\cong\ZZ$, and each divisorial ideal is represented by $T(a,0,0,0)$ where $a\in\ZZ$. 
In this case, there is a unique consistent dimer model associated with $R$. 
Since the existence of isoradial dimer model for each $3$-dimensional Gorenstein toric singularity is guaranteed by \cite{Gul}, 
such a unique consistent dimer model is isoradial. 
By the results in \cite[subsection~5.3]{Nak}, we have rank one MCM $R$-modules as in Figure~\ref{4b} from such an isoradial dimer model. 
In this figure, each circle stands for $a\in\Cl(R)$ corresponding to an MCM module $T(a,0,0,0)$. 
Especially, a double circle stands for the origin. 
Also, by the direct computation below, we can see all rank one MCM modules arising from an isoradial dimer model are conic. 

\begin{center}
\begin{tabular}{ll}
$T(1,0,0,0)\cong T(\lambda(\frac{1}{2},\frac{1}{2},0))$,& $T(2,0,0,0)\cong T(\lambda(\frac{1}{2},\frac{5}{4},\frac{1}{2})),$ \\ 
$T(3,0,0,0)\cong T(\lambda(\frac{3}{4},\frac{3}{2},\frac{5}{8}))$. 
\end{tabular} 
\end{center}

%%%%%%%%%%%%%%%%%%%%%%%%%%%%%%%%%%%%%%%%%%%%%%%%%%%%%%%
\subsubsection{{\bf \textup{Type 5a}}}
\label{type5a}

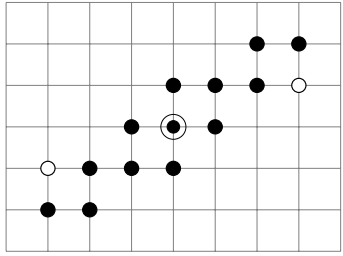
\begin{figure}[h]
{\scalebox{0.55}{
\begin{tikzpicture}
\draw [step=1,thin, gray] (-4,-3) grid (4,3);

\draw (0,0) circle [radius=0.3]; 
\filldraw (0,0) circle [radius=0.15]; 
\filldraw  [thick, fill=black] (0,1) circle [radius=0.17] ;\filldraw  [thick, fill=black] (1,0) circle [radius=0.17] ;
\filldraw  [thick, fill=black] (1,1) circle [radius=0.17] ;\filldraw  [thick, fill=black] (2,1) circle [radius=0.17] ;
\filldraw  [thick, fill=black] (2,2) circle [radius=0.17] ;\filldraw  [thick, fill=black] (3,2) circle [radius=0.17] ;
\filldraw  [thick, fill=black] (0,-1) circle [radius=0.17] ;\filldraw  [thick, fill=black] (-1,0) circle [radius=0.17] ;
\filldraw  [thick, fill=black] (-1,-1) circle [radius=0.17] ;\filldraw  [thick, fill=black] (-2,-1) circle [radius=0.17] ;
\filldraw  [thick, fill=black] (-2,-2) circle [radius=0.17] ;\filldraw  [thick, fill=black] (-3,-2) circle [radius=0.17] ;

\filldraw  [thick, fill=white] (3,1) circle [radius=0.17] ; \filldraw  [thick, fill=white] (-3,-1) circle [radius=0.17] ;

\end{tikzpicture}
}}
\caption{MCMs arising from consistent dimer models for Type 5a} 
\label{5a}
\end{figure}

We consider the reflexive polygon of type 5a. 
Thus, let $R$ be the $3$-dimensional complete local Gorenstein toric singularity defined by the cone $\sigma$: 
\[
\sigma=\mathrm{Cone}\{v_1=(1,0,1), v_2=(0,1,1), v_3=(-1,1,1), v_4=(-1,0,1), v_5=(0,-1,1) \}. 
\]
As elements in $\Cl(R)$, we obtain $2[I_1]+2[I_2]+[I_3]=0$, $3[I_1]+2[I_2]-[I_4]=0$, $2[I_1]+[I_2]+[I_5]=0$.
Thus, we have $\Cl(R)\cong\ZZ^2$, and each divisorial ideal is represented by $T(a,b,0,0)$ where $a, b\in\ZZ$. 
In this case, there are two consistent dimer models associated with $R$. 
By the results in \cite[subsection~5.4]{Nak}, we have rank one MCM $R$-modules as in Figure~\ref{5a} from such consistent dimer models. 
In this figure, each circle stands for $(a,b)\in\Cl(R)$ corresponding to an MCM module $T(a,b,0,0)$. 
Especially, a double circle stands for the origin $(0,0)$.
Also, by the direct computation below, we can see black circles are conic. 
On the other hand, white circles are not conic (e.g. check the condition in Lemma~\ref{conic_check}). 

\begin{center}
\begin{tabular}{ll}
$T(0,1,0,0)\cong T(\lambda(\frac{1}{8},\frac{3}{8},\frac{-1}{4}))$,& $T(1,0,0,0)\cong T(\lambda(\frac{1}{2},\frac{1}{4},\frac{-1}{4}))$,\\ 
$T(1,1,0,0)\cong T(\lambda(\frac{1}{2},\frac{1}{2},0))$,& $T(2,1,0,0)\cong T(\lambda(\frac{9}{8},\frac{1}{2},\frac{1}{4}))$, \\
$T(2,2,0,0)\cong T(\lambda(\frac{5}{4},\frac{3}{4},\frac{1}{2}))$,& $T(3,2,0,0)\cong T(\lambda(\frac{3}{2},\frac{3}{4},\frac{5}{8}))$.
\end{tabular} 
\end{center}

As we mentioned, there are two consistent dimer models associated with $R$, 
and one of them is isoradial and the other is not isoradial. 
The following figures are an isoradial one and the associated quiver. (The fraction on each arrow is an $\sfR$-charge which indicates this dimer model is isoradial.) 

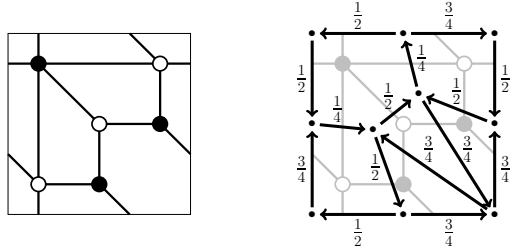
\begin{center}
\begin{tikzpicture} 
\node (DM) at (0,0) 
{\scalebox{0.4}{
\begin{tikzpicture}
%vertex
\node (W1) at (1,1){$$}; \node (B1) at (3,1){$$}; \node (W2) at (3,3){$$}; \node (B2) at (5,3){$$}; \node (W3) at (5,5){$$}; \node (B3) at (1,5){$$}; 
\draw[thick]  (0,0) rectangle (6,6);

%edge
\draw[line width=0.075cm]  (W1)--(B1); \draw[line width=0.075cm]  (W1)--(B3); \draw[line width=0.075cm]  (W1)--(1,0); 
\draw[line width=0.075cm]  (W1)--(0,2); \draw[line width=0.075cm]  (W2)--(B1); \draw[line width=0.075cm]  (W2)--(B2);
\draw[line width=0.075cm]  (W2)--(B3); \draw[line width=0.075cm]  (W3)--(B2); \draw[line width=0.075cm]  (W3)--(B3);
\draw[line width=0.075cm]  (W3)--(6,5); \draw[line width=0.075cm]  (W3)--(4,6); \draw[line width=0.075cm]  (4,0)--(B1); 
\draw[line width=0.075cm]  (6,2)--(B2); \draw[line width=0.075cm]  (0,5)--(B3); \draw[line width=0.075cm]  (1,6)--(B3);
%black
\filldraw  [ultra thick, fill=black] (3,1) circle [radius=0.24] ;\filldraw  [ultra thick, fill=black] (5,3) circle [radius=0.24] ;\filldraw  [ultra thick, fill=black] (1,5) circle [radius=0.24] ;
%white
\draw  [ultra thick,fill=white] (1,1) circle [radius=0.24] ;\draw  [ultra thick, fill=white] (3,3) circle [radius=0.24] ;\draw  [ultra thick, fill=white] (5,5) circle [radius=0.24] ;
\end{tikzpicture}
}}  ;

\node (Qui) at (4,0) 
{\scalebox{0.4}{ 
\begin{tikzpicture}
%vertex
\node (W1) at (1,1){$$}; \node (B1) at (3,1){$$}; \node (W2) at (3,3){$$}; \node (B2) at (5,3){$$};
\node (W3) at (5,5){$$}; \node (B3) at (1,5){$$}; 
\node (Q0) at (3.5,4){{\LARGE$\bullet$}}; \node (Q1) at (2,2.8){{\LARGE$\bullet$}}; \node (Q2a) at (0,3){{\LARGE$\bullet$}}; \node (Q2b) at (6,3){{\LARGE$\bullet$}};
\node (Q3a) at (6,0){{\LARGE$\bullet$}}; \node (Q3b) at (0,0){{\LARGE$\bullet$}}; \node (Q3c) at (0,6){{\LARGE$\bullet$}}; \node (Q3d) at (6,6){{\LARGE$\bullet$}}; 
\node (Q4a) at (3,0){{\LARGE$\bullet$}}; \node (Q4b) at (3,6){{\LARGE$\bullet$}};
%\node (Q0) at (3.5,4){{\LARGE$0$}}; \node (Q1) at (2,2.8){{\LARGE$1$}}; \node (Q2a) at (0,3){{\LARGE$2$}}; \node (Q2b) at (6,3){{\LARGE$2$}};
%\node (Q3a) at (6,0){{\LARGE$3$}}; \node (Q3b) at (0,0){{\LARGE$3$}}; \node (Q3c) at (0,6){{\LARGE$3$}}; \node (Q3d) at (6,6){{\LARGE$3$}}; 
%\node (Q4a) at (3,0){{\LARGE$4$}}; \node (Q4b) at (3,6){{\LARGE$4$}};
%\draw[thick]  (0,0) rectangle (6,6);
%edge
\draw[lightgray,line width=0.075cm]  (W1)--(B1); \draw[lightgray,line width=0.075cm]  (W1)--(B3); \draw[lightgray,line width=0.075cm]  (W1)--(1,0); 
\draw[lightgray,line width=0.075cm]  (W1)--(0,2); \draw[lightgray,line width=0.075cm]  (W2)--(B1); \draw[lightgray,line width=0.075cm]  (W2)--(B2);
\draw[lightgray,line width=0.075cm]  (W2)--(B3); \draw[lightgray,line width=0.075cm]  (W3)--(B2); \draw[lightgray,line width=0.075cm]  (W3)--(B3);
\draw[lightgray,line width=0.075cm]  (W3)--(6,5); \draw[lightgray,line width=0.075cm]  (W3)--(4,6); \draw[lightgray,line width=0.075cm]  (4,0)--(B1); 
\draw[lightgray,line width=0.075cm]  (6,2)--(B2); \draw[lightgray,line width=0.075cm]  (0,5)--(B3); \draw[lightgray,line width=0.075cm]  (1,6)--(B3);

%black
\filldraw  [ultra thick, draw=lightgray, fill=lightgray] (3,1) circle [radius=0.24] ;\filldraw  [ultra thick, draw=lightgray, fill=lightgray] (5,3) circle [radius=0.24] ;
\filldraw  [ultra thick, draw=lightgray, fill=lightgray] (1,5) circle [radius=0.24] ;
%white
\draw  [ultra thick, draw=lightgray, fill=white] (1,1) circle [radius=0.24] ;\draw  [ultra thick, draw=lightgray, fill=white] (3,3) circle [radius=0.24] ;
\draw  [ultra thick, draw=lightgray, fill=white] (5,5) circle [radius=0.24] ;

\draw[->, line width=0.1cm] (Q1)--node[midway,xshift=-7pt,yshift=15pt] {{\Huge$\frac{1}{2}$}}(Q0); 
\draw[->, line width=0.1cm] (Q0)--node[midway,xshift=10pt,yshift=5pt] {{\Huge$\frac{3}{4}$}}(Q3a); 
\draw[->, line width=0.1cm] (Q3a)--node[midway,xshift=-5pt,yshift=22pt] {{\Huge$\frac{3}{4}$}}(Q1); 
\draw[->, line width=0.1cm] (Q4a)--node[midway,xshift=0pt,yshift=-17pt] {{\Huge$\frac{3}{4}$}}(Q3a);
\draw[->, line width=0.1cm] (Q1)--node[midway,xshift=-12pt,yshift=5pt] {{\Huge$\frac{1}{2}$}}(Q4a); 
\draw[->, line width=0.1cm] (Q3a)--node[midway,xshift=10pt,yshift=0pt] {{\Huge$\frac{3}{4}$}}(Q2b); 
\draw[->, line width=0.1cm] (Q2b)--node[midway,xshift=0pt,yshift=17pt] {{\Huge$\frac{1}{2}$}}(Q0);
\draw[->, line width=0.1cm] (Q0)--node[midway,xshift=10pt,yshift=5pt] {{\Huge$\frac{1}{4}$}}(Q4b); 
\draw[->, line width=0.1cm] (Q2a)--node[midway,xshift=-5pt,yshift=17pt] {{\Huge$\frac{1}{4}$}}(Q1); 
\draw[->, line width=0.1cm] (Q4a)--node[midway,xshift=0pt,yshift=-17pt] {{\Huge$\frac{1}{2}$}}(Q3b);
\draw[->, line width=0.1cm] (Q3b)--node[midway,xshift=-10pt,yshift=0pt] {{\Huge$\frac{3}{4}$}}(Q2a);  
\draw[->, line width=0.1cm] (Q3c)--node[midway,xshift=-10pt,yshift=0pt] {{\Huge$\frac{1}{2}$}}(Q2a); 
\draw[->, line width=0.1cm] (Q4b)--node[midway,xshift=0pt,yshift=17pt] {{\Huge$\frac{1}{2}$}}(Q3c); 
\draw[->, line width=0.1cm] (Q4b)--node[midway,xshift=0pt,yshift=17pt] {{\Huge$\frac{3}{4}$}}(Q3d);
\draw[->, line width=0.1cm] (Q3d)--node[midway,xshift=10pt,yshift=0pt] {{\Huge$\frac{1}{2}$}}(Q2b);  
\end{tikzpicture}
}} ;

\end{tikzpicture}
\end{center}
Again, by the results in \cite[subsection~5.4]{Nak}, we can see rank one MCM modules arising from this isoradial dimer model are actually conic, 
and non-conic divisorial ideals arise from a consistent dimer model which is not isoradial. 

%%%%%%%%%%%%%%%%%%%%%%%%%%%%%%%%%%%%%%%%%%%%%%%%%%%%%%%
\subsubsection{{\bf \textup{Type 5b}}}
\label{type5b}

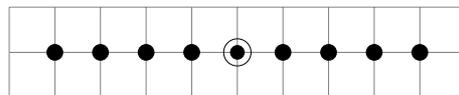
\begin{figure}[h]
{\scalebox{0.6}{
\begin{tikzpicture}
\draw [step=1,thin, gray] (-5,-1) grid (5,1);

\draw (0,0) circle [radius=0.3]; 
\filldraw (0,0) circle [radius=0.15]; 
\filldraw  [thick, fill=black] (1,0) circle [radius=0.17] ;
\filldraw  [thick, fill=black] (2,0) circle [radius=0.17] ;
\filldraw  [thick, fill=black] (3,0) circle [radius=0.17] ;
\filldraw  [thick, fill=black] (4,0) circle [radius=0.17] ;
\filldraw  [thick, fill=black] (-1,0) circle [radius=0.17] ;
\filldraw  [thick, fill=black] (-2,0) circle [radius=0.17] ;
\filldraw  [thick, fill=black] (-3,0) circle [radius=0.17] ;
\filldraw  [thick, fill=black] (-4,0) circle [radius=0.17] ;
%\node  at (1,-1.8) {{\Large\text{MCMs arising from consistent dimer models}}} ;
\end{tikzpicture}
}}
\caption{MCMs arising from consistent dimer models for Type 5b} 
\label{5b}
\end{figure}

We consider the reflexive polygon of type 5b. 
Thus, let $R$ be the $3$-dimensional complete local Gorenstein toric singularity defined by the cone $\sigma$: 
\[
\sigma=\mathrm{Cone}\{v_1=(0,1,1), v_2=(-1,0,1), v_3=(-1,-1,1), v_4=(1,-1,1) \}. 
\]
As elements in $\Cl(R)$, we obtain $[I_1]+2[I_4]=0$, $[I_2]-4[I_4]=0$, $[I_3]+3[I_4]=0$.
Thus, we have $\Cl(R)\cong\ZZ$, and each divisorial ideal is represented by $T(0,0,0,d)$ where $d\in\ZZ$. 
In this case, there is a unique consistent dimer model associated with $R$.  
Since the existence of isoradial dimer model for each $3$-dimensional Gorenstein toric singularity is guaranteed by \cite{Gul}, 
such a unique consistent dimer model is isoradial. 
By the results in \cite[subsection~5.5]{Nak}, we have rank one MCM $R$-modules as in Figure~\ref{5b} from such an isoradial dimer model. 
In this figure, each circle stands for $d\in\Cl(R)$ corresponding to an MCM module $T(0,0,0,d)$. 
Especially, a double circle stands for the origin. 
Also, by the direct computation below, we can see all rank one MCM modules arising from an isoradial dimer model are conic. 

\begin{center}
\begin{tabular}{lll}
$T(0,0,0,1)\cong T(\lambda(\frac{1}{4},\frac{-3}{8},\frac{-1}{2}))$,& $T(0,0,0,2)\cong T(\lambda(\frac{3}{4},\frac{-5}{8},\frac{-1}{8}))$,\\
$T(0,0,0,3)\cong T(\lambda(\frac{5}{4},\frac{-1}{2},\frac{1}{2}))$,& $T(0,0,0,4)\cong T(\lambda(\frac{13}{8},\frac{-7}{8},\frac{3}{4}))$.  
\end{tabular} 
\end{center}

%%%%%%%%%%%%%%%%%%%%%%%%%%%%%%%%%%%%%%%%%%%%%%%%%%%%%%%
\subsubsection{{\bf \textup{Type 6a}}}
\label{type6a}

We consider the reflexive polygon of type 6a. 
Thus, let $R$ be the $3$-dimensional complete local Gorenstein toric singularity defined by the cone $\sigma$: 
\[
\sigma=\mathrm{Cone}\{v_1=(1,0,1), v_2=(0,1,1), v_3=(-1,1,1), v_4=(-1,0,1), v_5=(0,-1,1), v_6=(1,-1,1) \}. 
\]
As an element in $\Cl(R)$, we obtain $[I_1]+2[I_2]+2[I_3]+[I_4]=0$, $2[I_1]+3[I_2]+2[I_3]-[I_5]=0$, $2[I_1]+2[I_2]+[I_3]+[I_6]=0$.  
Therefore, we have $\Cl(R)\cong\ZZ^3$, and each divisorial ideal is represented by $T(a,b,c,0,0,0)$ where $a, b, c\in\ZZ$. 
In this case, there are five consistent dimer models associated with $R$. 
By the results in \cite[subsection~5.6]{Nak}, rank one MCM $R$-modules arising from such consistent dimer models are 
modules corresponding to elements $(a,b,c)=(0,0,0)$, $(1,0,0)$, $(0,1,0)$, $(0,0,1)$, $(1,1,0)$, $(0,1,1)$, $(1,1,1)$, $(1,2,1)$, $(1,2,2)$, $(2,2,1)$, $(2,3,2)$, 
$(2,2,2)$, $(2,3,1)$, $(1,0,-1)$, $(1,3,1)$ in $\Cl(R)$ and their $R$-duals. 
Also, we can see MCM modules listed below and their $R$-duals are non-free conic divisorial ideals. 

\medskip

\begin{center}
\begin{tabular}{ll}
$T(1,0,0,0,0,0)\cong T(\lambda(\frac{1}{2},\frac{1}{4},\frac{-3}{8}))$,& $T(0,1,0,0,0,0)\cong T(\lambda(\frac{1}{4},\frac{1}{2},\frac{-3}{8}))$, \\
$T(0,0,1,0,0,0)\cong T(\lambda(\frac{-1}{4},\frac{1}{4},\frac{-3}{8}))$,& $T(1,1,0,0,0,0)\cong T(\lambda(\frac{1}{2},\frac{1}{2},\frac{-1}{4}))$, \\
$T(0,1,1,0,0,0)\cong T(\lambda(0,\frac{1}{2},\frac{-1}{4}))$,& $T(1,1,1,0,0,0)\cong T(\lambda(\frac{3}{8},\frac{11}{16},\frac{-1}{4}))$, \\
$T(1,2,1,0,0,0)\cong T(\lambda(\frac{1}{2},1,\frac{1}{4}))$,& $T(1,2,2,0,0,0)\cong T(\lambda(\frac{1}{2},\frac{5}{4},\frac{3}{8}))$, \\
$T(2,2,1,0,0,0)\cong T(\lambda(\frac{3}{4},\frac{5}{4},\frac{3}{8}))$,& $T(2,3,2,0,0,0)\cong T(\lambda(\frac{3}{4},\frac{3}{2},\frac{5}{8}))$, \\
$T(2,2,2,0,0,0)\cong T(\lambda(\frac{5}{8},\frac{11}{8},\frac{1}{2}))$.
\end{tabular} 
\end{center}

On the other hand, we can check remaining ones are not conic (e.g. check the condition in Lemma~\ref{conic_check}). 

As we mentioned, there are five consistent dimer models associated with $R$. One of them is isoradial and the others are not isoradial. 
The following figures are an isoradial one and the associated quiver. 
(The fraction on each arrow is an $\sfR$-charge which indicates this dimer model is isoradial.)

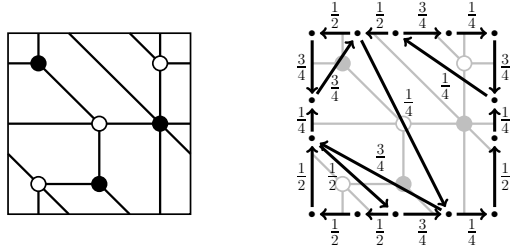
\begin{center}
\begin{tikzpicture} 
\node (DM) at (0,0) 
{\scalebox{0.4}{
\begin{tikzpicture}
%vertex
\node (B1) at (3,1){$$}; \node (B2) at (5,3){$$}; \node (B3) at (1,5){$$};  
\node (W1) at (1,1){$$}; \node (W2) at (3,3){$$}; \node (W3) at (5,5){$$}; 
\draw[ultra thick]  (0,0) rectangle (6,6);
%edge
\draw[line width=0.075cm]  (W1)--(1,0); \draw[line width=0.075cm]  (W1)--(2,0); \draw[line width=0.075cm]  (W1)--(0,2);
\draw[line width=0.075cm]  (B1)--(W1); \draw[line width=0.075cm]  (B1)--(W2); \draw[line width=0.075cm]  (B1)--(4,0);
\draw[line width=0.075cm]  (W2)--(0,3); \draw[line width=0.075cm]  (W2)--(B2); \draw[line width=0.075cm]  (W2)--(B3);
\draw[line width=0.075cm]  (B2)--(5,0); \draw[line width=0.075cm]  (B2)--(6,2); \draw[line width=0.075cm]  (B2)--(6,3); 
\draw[line width=0.075cm]  (B2)--(W3); \draw[line width=0.075cm]  (B2)--(2,6); \draw[line width=0.075cm]  (B3)--(0,5); 
\draw[line width=0.075cm]  (B3)--(1,6); \draw[line width=0.075cm]  (W3)--(4,6); \draw[line width=0.075cm]  (W3)--(5,6); 
\draw[line width=0.075cm]  (W3)--(6,5);
%black
\filldraw  [ultra thick, fill=black] (3,1) circle [radius=0.24] ;\filldraw  [ultra thick, fill=black] (5,3) circle [radius=0.24] ;
\filldraw  [ultra thick, fill=black] (1,5) circle [radius=0.24] ;
%white
\draw  [ultra thick,fill=white] (1,1) circle [radius=0.24] ;\draw  [ultra thick, fill=white] (3,3) circle [radius=0.24] ;
\draw  [ultra thick,fill=white] (5,5)circle [radius=0.24] ; 
\end{tikzpicture}
}}  ;

\node (Qui) at (4,0) 
{\scalebox{0.4}{ 
\begin{tikzpicture}
%vertex
\node (B1) at (3,1){$$}; \node (B2) at (5,3){$$}; \node (B3) at (1,5){$$};  
\node (W1) at (1,1){$$}; \node (W2) at (3,3){$$}; \node (W3) at (5,5){$$}; 
\node (Q0a) at (2.75,0){{\LARGE$\bullet$}}; \node (Q0b) at (2.75,6){{\LARGE$\bullet$}}; \node (Q1a) at (1.5,0){{\LARGE$\bullet$}}; 
\node (Q1b) at (1.5,6){{\LARGE$\bullet$}}; 
\node (Q2a) at (0,0){{\LARGE$\bullet$}}; \node (Q2b) at (6,0){{\LARGE$\bullet$}}; \node (Q2c) at (6,6){{\LARGE$\bullet$}}; \node (Q2d) at (0,6){{\LARGE$\bullet$}}; 
\node (Q3a) at (0,3.75){{\LARGE$\bullet$}}; \node (Q3b) at (6,3.75){{\LARGE$\bullet$}};
\node (Q4a) at (4.5,0){{\LARGE$\bullet$}}; \node (Q4b) at (4.5,6){{\LARGE$\bullet$}}; \node (Q5a) at (0,2.5){{\LARGE$\bullet$}}; 
\node (Q5b) at (6,2.5){{\LARGE$\bullet$}};

%\node (Q0a) at (2.75,0){{\LARGE$0$}}; \node (Q0b) at (2.75,6){{\LARGE$0$}}; \node (Q1a) at (1.5,0){{\LARGE$1$}}; \node (Q1b) at (1.5,6){{\LARGE$1$}}; 
%\node (Q2a) at (0,0){{\LARGE$2$}}; \node (Q2b) at (6,0){{\LARGE$2$}}; \node (Q2c) at (6,6){{\LARGE$2$}}; \node (Q2d) at (0,6){{\LARGE$2$}}; 
%\node (Q3a) at (0,3.75){{\LARGE$3$}}; \node (Q3b) at (6,3.75){{\LARGE$3$}};
%\node (Q4a) at (4.5,0){{\LARGE$4$}}; \node (Q4b) at (4.5,6){{\LARGE$4$}}; \node (Q5a) at (0,2.5){{\LARGE$5$}}; \node (Q5b) at (6,2.5){{\LARGE$5$}};
%\draw[ultra thick]  (0,0) rectangle (6,6);

%edge
\draw[lightgray, line width=0.075cm]  (W1)--(1,0); \draw[lightgray, line width=0.075cm]  (W1)--(2,0); \draw[lightgray, line width=0.075cm]  (W1)--(0,2);
\draw[lightgray, line width=0.075cm]  (B1)--(W1); \draw[lightgray, line width=0.075cm]  (B1)--(W2); \draw[lightgray, line width=0.075cm]  (B1)--(4,0);
\draw[lightgray, line width=0.075cm]  (W2)--(0,3); \draw[lightgray, line width=0.075cm]  (W2)--(B2); \draw[lightgray, line width=0.075cm]  (W2)--(B3);
\draw[lightgray, line width=0.075cm]  (B2)--(5,0); \draw[lightgray, line width=0.075cm]  (B2)--(6,2); \draw[lightgray, line width=0.075cm]  (B2)--(6,3); 
\draw[lightgray, line width=0.075cm]  (B2)--(W3); \draw[lightgray, line width=0.075cm]  (B2)--(2,6); \draw[lightgray, line width=0.075cm]  (B3)--(0,5); 
\draw[lightgray, line width=0.075cm]  (B3)--(1,6); \draw[lightgray, line width=0.075cm]  (W3)--(4,6); \draw[lightgray, line width=0.075cm]  (W3)--(5,6); 
\draw[lightgray, line width=0.075cm]  (W3)--(6,5);
%black
\filldraw  [ultra thick, draw=lightgray, fill=lightgray] (3,1) circle [radius=0.24] ;\filldraw  [ultra thick, draw=lightgray, fill=lightgray] (5,3) circle [radius=0.24] ;
\filldraw  [ultra thick, draw=lightgray, fill=lightgray] (1,5) circle [radius=0.24] ;
%white
\draw  [ultra thick,draw=lightgray,fill=white] (1,1) circle [radius=0.24] ;\draw  [ultra thick, draw=lightgray,fill=white] (3,3) circle [radius=0.24] ;
\draw  [ultra thick,draw=lightgray,fill=white] (5,5)circle [radius=0.24] ; 
%arrow
\draw[->, line width=0.1cm] (Q0b)--node[midway,xshift=3pt,yshift=17pt] {{\Huge$\frac{1}{2}$}}(Q1b); 
\draw[->, line width=0.1cm] (Q0a)--node[midway,xshift=0pt,yshift=-17pt] {{\Huge$\frac{3}{4}$}}(Q4a); 
\draw[->, line width=0.1cm] (Q0b)--node[midway,xshift=0pt,yshift=17pt] {{\Huge$\frac{3}{4}$}}(Q4b); 
\draw[->, line width=0.1cm] (Q4b)--node[midway,xshift=0pt,yshift=17pt] {{\Huge$\frac{1}{4}$}}(Q2c); 
\draw[->, line width=0.1cm] (Q1b)--node[midway,xshift=0pt,yshift=17pt] {{\Huge$\frac{1}{2}$}}(Q2d); 
\draw[->, line width=0.1cm] (Q4a)--node[midway,xshift=0pt,yshift=-17pt] {{\Huge$\frac{1}{4}$}}(Q2b); 
\draw[->, line width=0.1cm] (Q1b)--node[midway,xshift=5pt,yshift=20pt] {{\Huge$\frac{1}{4}$}}(Q4a);
\draw[->, line width=0.1cm] (Q2a)--node[midway,xshift=-10pt,yshift=0pt] {{\Huge$\frac{1}{2}$}}(Q5a); 
\draw[->, line width=0.1cm] (Q2c)--node[midway,xshift=10pt,yshift=0pt] {{\Huge$\frac{3}{4}$}}(Q3b); 
\draw[->, line width=0.1cm] (Q2d)--node[midway,xshift=-10pt,yshift=0pt] {{\Huge$\frac{3}{4}$}}(Q3a);
\draw[->, line width=0.1cm] (Q3a)--node[midway,xshift=0pt,yshift=-20pt] {{\Huge$\frac{3}{4}$}}(Q1b); 
\draw[->, line width=0.1cm] (Q3b)--node[midway,xshift=0pt,yshift=-18pt] {{\Huge$\frac{1}{4}$}}(Q0b); 
\draw[->, line width=0.1cm] (Q4a)--node[midway,xshift=0pt,yshift=18pt] {{\Huge$\frac{3}{4}$}}(Q5a);
\draw[->, line width=0.1cm] (Q5a)--node[midway,xshift=-20pt,yshift=0pt] {{\Huge$\frac{1}{2}$}}(Q0a); 
\draw[->, line width=0.1cm] (Q5a)--node[midway,xshift=-10pt,yshift=0pt] {{\Huge$\frac{1}{4}$}}(Q3a); 
\draw[->, line width=0.1cm] (Q2b)--node[midway,xshift=10pt,yshift=0pt] {{\Huge$\frac{1}{2}$}}(Q5b); 
\draw[->, line width=0.1cm] (Q5b)--node[midway,xshift=10pt,yshift=0pt] {{\Huge$\frac{1}{4}$}}(Q3b); 
\draw[->, line width=0.1cm] (Q0a)--node[midway,xshift=3pt,yshift=-17pt] {{\Huge$\frac{1}{2}$}}(Q1a); 
\draw[->, line width=0.1cm] (Q1a)--node[midway,xshift=0pt,yshift=-17pt] {{\Huge$\frac{1}{2}$}}(Q2a);
\end{tikzpicture}
}} ;

\end{tikzpicture}
\end{center}
Again, by the results in \cite[subsection~5.6]{Nak}, we can see rank one MCM modules arising from this isoradial dimer model are actually conic, 
and non-conic divisorial ideals arise from consistent dimer models which are not isoradial. 
Further, we remark that even although $T(2,2,2,0,0,0)$ and its $R$-dual are conic, they do not arise from an isoradial dimer model. 
They arise from other consistent dimer models.

%%%%%%%%%%%%%%%%%%%%%%%%%%%%%%%%%%%%%%%%%%%%%%%%%%%%%%%
\subsubsection{{\bf \textup{Type 6b}}}
\label{type6b}

We consider the reflexive polygon of type 6b. 
Thus, let $R$ be the $3$-dimensional complete local Gorenstein toric singularity defined by the cone $\sigma$: 
\[
\sigma=\mathrm{Cone}\{v_1=(1,0,1), v_2=(0,1,1), v_3=(-1,1,1), v_4=(-1,-1,1), v_5=(0,-1,1) \}. 
\]
As elements in $\Cl(R)$, we obtain $[I_1]+2[I_4]+2[I_5]=0$, $[I_2]-4[I_4]-3[I_5]=0$, $[I_3]+3[I_4]+2[I_5]=0$. 
Thus, we have $\Cl(R)\cong\ZZ^2$, and each divisorial ideal is represented by $T(0,0,0,d,e)$ where $d, e\in\ZZ$. 
In this case, there are three consistent dimer models associated with $R$. 
By the results in \cite[subsection~5.7]{Nak}, we have rank one MCM $R$-modules as in Figure~\ref{6b} from such consistent dimer models. 
In this figure, each circle stands for $(d,e)\in\Cl(R)$ corresponding to an MCM module $T(0,0,0,d,e)$. 
Especially, a double circle stands for the origin $(0,0)$. Also, by the direct computation below, we can see black circles are conic. 
On the other hand, white circles are not conic (e.g. check the condition in Lemma~\ref{conic_check}). 

\begin{center}
\begin{tabular}{ll}
$T(0,0,0,0,1)\cong T(\lambda(\frac{1}{4},\frac{-5}{16},\frac{-1}{4}))$,& $T(0,0,0,1,0)\cong T(\lambda(\frac{-1}{4},\frac{-1}{8},\frac{-1}{4}))$,\\
$T(0,0,0,1,1)\cong T(\lambda(\frac{-1}{4},\frac{-1}{2},\frac{-1}{4}))$,& $T(0,0,0,2,1)\cong T(\lambda(\frac{-11}{16},\frac{-11}{16},\frac{-1}{4}))$, \\
$T(0,0,0,2,2)\cong T(\lambda(\frac{-3}{8},\frac{-7}{8},\frac{1}{4}))$,& $T(0,0,0,3,2)\cong T(\lambda(\frac{-3}{4},\frac{-9}{8},\frac{1}{4}))$, \\
$T(0,0,0,3,3)\cong T(\lambda(\frac{-3}{4},\frac{-3}{2},\frac{5}{8}))$,& $T(0,0,0,4,3)\cong T(\lambda(\frac{-13}{16},\frac{-27}{16},\frac{3}{4}))$.
\end{tabular} 
\end{center}

\begin{figure}[h]
{\scalebox{0.55}{
\begin{tikzpicture}
\draw [step=1,thin, gray] (-5,-4) grid (5,4);

\draw (0,0) circle [radius=0.3]; 
\filldraw (0,0) circle [radius=0.15]; 
\filldraw  [thick, fill=black] (0,1) circle [radius=0.17] ;\filldraw  [thick, fill=black] (1,0) circle [radius=0.17] ;
\filldraw  [thick, fill=black] (1,1) circle [radius=0.17] ;\filldraw  [thick, fill=black] (2,1) circle [radius=0.17] ;
\filldraw  [thick, fill=black] (2,2) circle [radius=0.17] ;\filldraw  [thick, fill=black] (3,2) circle [radius=0.17] ;\filldraw  [thick, fill=black] (4,3) circle [radius=0.17] ;
\filldraw  [thick, fill=black] (0,-1) circle [radius=0.17] ;\filldraw  [thick, fill=black] (-1,0) circle [radius=0.17] ;
\filldraw  [thick, fill=black] (-1,-1) circle [radius=0.17] ;\filldraw  [thick, fill=black] (-2,-1) circle [radius=0.17] ;
\filldraw  [thick, fill=black] (-2,-2) circle [radius=0.17] ;\filldraw  [thick, fill=black] (-3,-2) circle [radius=0.17] ;\filldraw  [thick, fill=black] (-4,-3) circle [radius=0.17] ;
\filldraw  [thick, fill=black] (3,3) circle [radius=0.17] ;\filldraw  [thick, fill=black] (-3,-3) circle [radius=0.17] ;

\filldraw  [thick, fill=white] (4,2) circle [radius=0.17] ; \filldraw  [thick, fill=white] (-4,-2) circle [radius=0.17] ;
\end{tikzpicture}
}}
\caption{MCMs arising from consistent dimer models for Type 6b} 
\label{6b}
\end{figure}

As we mentioned, there are three consistent dimer models associated with $R$, 
and two of them are isoradial and the other is not isoradial. 
The following figures are isoradial ones and the associated quivers. (The fraction on each arrow is an $\sfR$-charge which indicates they are isoradial.)

\begin{center}
\begin{tikzpicture} 
\node (DM1) at (0,0) 
{\scalebox{0.6}{
\begin{tikzpicture}
%vertex
\node (B1) at (0.5,0.5){$$}; \node (B2) at (2.5,1.5){$$}; \node (B3) at (0.5,2.5){$$};  \node (B4) at (2.5,3.5){$$};
\node (W1) at (3.5,0.5){$$}; \node (W2) at (1.5,1.5){$$}; \node (W3) at (3.5,2.5){$$};  \node (W4) at (1.5,3.5){$$};
\draw[thick]  (0,0) rectangle (4,4);

%edge
\draw[line width=0.05cm]  (B1)--(W2); \draw[line width=0.05cm]  (B1)--(0,0.5); \draw[line width=0.05cm]  (B1)--(1,0); 
\draw[line width=0.05cm]  (W2)--(B2); \draw[line width=0.05cm]  (W2)--(B3); \draw[line width=0.05cm]  (W3)--(B2);
\draw[line width=0.05cm]  (B2)--(W1); \draw[line width=0.05cm]  (W4)--(B2); \draw[line width=0.05cm]  (B2)--(1.75,0); 
\draw[line width=0.05cm]  (W1)--(3,0); \draw[line width=0.05cm]  (W1)--(4,0.5); \draw[line width=0.05cm]  (B3)--(W4); 
\draw[line width=0.05cm]  (B3)--(0,2.5); \draw[line width=0.05cm]  (W3)--(4,2.5); \draw[line width=0.05cm]  (W3)--(B4);
\draw[line width=0.05cm]  (W4)--(B4); \draw[line width=0.05cm]  (W4)--(1,4); \draw[line width=0.05cm]  (W4)--(1.75,4); \draw[line width=0.05cm]  (B4)--(3,4);
%black
\filldraw  [ultra thick, fill=black] (0.5,0.5) circle [radius=0.16] ;
\filldraw  [ultra thick, fill=black] (2.5,1.5) circle [radius=0.16] ;
\filldraw  [ultra thick, fill=black] (0.5,2.5) circle [radius=0.16] ;
\filldraw  [ultra thick, fill=black] (2.5,3.5) circle [radius=0.16] ;
%white
\draw  [ultra thick,fill=white] (3.5,0.5) circle [radius=0.16] ;
\draw  [ultra thick, fill=white] (1.5,1.5) circle [radius=0.16] ;
\draw  [ultra thick,fill=white] (3.5,2.5)circle [radius=0.16] ;
\draw  [ultra thick, fill=white] (1.5,3.5) circle [radius=0.16] ;
\end{tikzpicture}
}}  ;

\node (Qui1) at (4,0) 
{\scalebox{0.6}{ 
\begin{tikzpicture}
%vertex
\node (B1) at (0.5,0.5){$$}; \node (B2) at (2.5,1.5){$$}; \node (B3) at (0.5,2.5){$$};  \node (B4) at (2.5,3.5){$$};
\node (W1) at (3.5,0.5){$$}; \node (W2) at (1.5,1.5){$$}; \node (W3) at (3.5,2.5){$$};  \node (W4) at (1.5,3.5){$$};
%\draw[thick]  (0,0) rectangle (4,4);
\node (Q0a) at (0,1.5){$\bullet$}; \node (Q0b) at (4,1.5){$\bullet$}; 
\node (Q1) at (1.5,2.2){$\bullet$}; \node (Q2a) at (1.25,0){$\bullet$}; \node (Q2b) at (1.25,4){$\bullet$}; 
\node (Q3a) at (2.5,0){$\bullet$}; \node (Q3b) at (2.5,4){$\bullet$}; \node (Q4) at (2.6,2.5){$\bullet$}; 
\node (Q5a) at (0,0){$\bullet$}; \node (Q5b) at (4,0){$\bullet$}; \node (Q5c) at (4,4){$\bullet$}; \node (Q5d) at (0,4){$\bullet$}; 

%\node (Q0a) at (0,1.5){$0$}; \node (Q0b) at (4,1.5){$0$}; 
%\node (Q1) at (1.5,2.2){$1$}; \node (Q2a) at (1.25,0){$2$}; \node (Q2b) at (1.25,4){$2$}; 
%\node (Q3a) at (2.5,0){$3$}; \node (Q3b) at (2.5,4){$3$}; \node (Q4) at (2.6,2.5){$4$}; 
%\node (Q5a) at (0,0){$5$}; \node (Q5b) at (4,0){$5$}; \node (Q5c) at (4,4){$5$}; \node (Q5d) at (0,4){$5$}; 
%edge
\draw[lightgray, line width=0.05cm]  (B1)--(W2); \draw[lightgray, line width=0.05cm]  (B1)--(0,0.5); \draw[lightgray, line width=0.05cm]  (B1)--(1,0); 
\draw[lightgray, line width=0.05cm]  (W2)--(B2); \draw[lightgray, line width=0.05cm]  (W2)--(B3); \draw[lightgray, line width=0.05cm]  (W3)--(B2);
\draw[lightgray, line width=0.05cm]  (B2)--(W1); \draw[lightgray, line width=0.05cm]  (W4)--(B2); \draw[lightgray, line width=0.05cm]  (B2)--(1.75,0); 
\draw[lightgray, line width=0.05cm]  (W1)--(3,0); \draw[lightgray, line width=0.05cm]  (W1)--(4,0.5); \draw[lightgray, line width=0.05cm]  (B3)--(W4); 
\draw[lightgray, line width=0.05cm]  (B3)--(0,2.5); \draw[lightgray, line width=0.05cm]  (W3)--(4,2.5); \draw[lightgray, line width=0.05cm]  (W3)--(B4);
\draw[lightgray, line width=0.05cm]  (W4)--(B4); \draw[lightgray, line width=0.05cm]  (W4)--(1,4); \draw[lightgray, line width=0.05cm]  (W4)--(1.75,4); 
\draw[lightgray, line width=0.05cm]  (B4)--(3,4);
%black
\filldraw  [ultra thick, draw=lightgray, fill=lightgray] (0.5,0.5) circle [radius=0.16] ;\filldraw  [ultra thick, draw=lightgray, fill=lightgray] (2.5,1.5) circle [radius=0.16] ;
\filldraw  [ultra thick, draw=lightgray, fill=lightgray] (0.5,2.5) circle [radius=0.16] ;\filldraw  [ultra thick, draw=lightgray, fill=lightgray] (2.5,3.5) circle [radius=0.16] ;
%white
\draw  [ultra thick,draw=lightgray,fill=white] (3.5,0.5) circle [radius=0.16] ;\draw  [ultra thick, draw=lightgray,fill=white] (1.5,1.5) circle [radius=0.16] ;
\draw  [ultra thick,draw=lightgray,fill=white] (3.5,2.5)circle [radius=0.16] ;\draw  [ultra thick, draw=lightgray,fill=white] (1.5,3.5) circle [radius=0.16] ;

\draw[->, line width=0.067cm] (Q5a)--node[midway,xshift=0pt,yshift=-15pt] {{\LARGE$\frac{1}{2}$}}(Q2a); 
\draw[->, line width=0.067cm] (Q5b)--node[midway,xshift=0pt,yshift=-15pt] {{\LARGE$\frac{3}{4}$}}(Q3a); 
\draw[->, line width=0.067cm] (Q5d)--node[midway,xshift=0pt,yshift=15pt] {{\LARGE$\frac{1}{2}$}}(Q2b); 
\draw[->, line width=0.067cm] (Q5c)--node[midway,xshift=0pt,yshift=15pt] {{\LARGE$\frac{3}{4}$}}(Q3b); 
\draw[->, line width=0.067cm] (Q5c)--node[midway,xshift=7pt,yshift=0pt] {{\LARGE$\frac{3}{4}$}}(Q0b); 
\draw[->, line width=0.067cm] (Q5d)--node[midway,xshift=-7pt,yshift=0pt] {{\LARGE$\frac{3}{4}$}}(Q0a);   
\draw[->, line width=0.067cm] (Q1)--node[midway,xshift=7pt,yshift=5pt] {{\LARGE$\frac{1}{2}$}}(Q5d); 
\draw[->, line width=0.067cm] (Q1)--node[midway,xshift=7pt,yshift=0pt] {{\LARGE$\frac{1}{2}$}}(Q2a); 
\draw[->, line width=0.067cm] (Q2a)--node[midway,xshift=0pt,yshift=-15pt] {{\LARGE$\frac{1}{4}$}}(Q3a); 
\draw[->, line width=0.067cm] (Q2a)--node[midway,xshift=3pt,yshift=10pt] {{\LARGE$\frac{3}{4}$}}(Q0a); 
\draw[->, line width=0.067cm] (Q2b)--node[midway,xshift=0pt,yshift=15pt] {{\LARGE$\frac{1}{4}$}}(Q3b); 
\draw[->, line width=0.067cm] (Q3b)--node[midway,xshift=-7pt,yshift=0pt] {{\LARGE$\frac{1}{2}$}}(Q4); 
\draw[->, line width=0.067cm] (Q3a)--node[midway,xshift=-3pt,yshift=10pt] {{\LARGE$\frac{1}{2}$}}(Q0b); 
\draw[->, line width=0.067cm] (Q4)--node[midway,xshift=-7pt,yshift=5pt] {{\LARGE$\frac{3}{4}$}}(Q5c); 
\draw[->, line width=0.067cm] (Q4)--node[midway,xshift=3pt,yshift=-12pt] {{\LARGE$\frac{1}{4}$}}(Q1);  
\draw[->, line width=0.067cm] (Q0a)--node[midway,xshift=-7pt,yshift=3pt] {{\LARGE$\frac{3}{4}$}}(Q5a); 
\draw[->, line width=0.067cm] (Q0b)--node[midway,xshift=7pt,yshift=3pt] {{\LARGE$\frac{3}{4}$}}(Q5b); 
\draw[->, line width=0.067cm] (Q0a)--node[midway,xshift=0pt,yshift=12pt] {{\LARGE$\frac{3}{4}$}}(Q1); 
\draw[->, line width=0.067cm] (Q0b)--node[midway,xshift=3pt,yshift=10pt] {{\LARGE$\frac{1}{2}$}}(Q4);  
\end{tikzpicture}
}} ;

\end{tikzpicture}
\end{center}

%\medskip

\begin{center}
\begin{tikzpicture} 
\node (DM2) at (0,0) 
{\scalebox{0.6}{
\begin{tikzpicture}
%vertex
\node (B1) at (0.5,0.5){$$}; \node (B2) at (2.5,1.5){$$}; \node (B3) at (0.5,2.5){$$};  \node (B4) at (2.5,3.5){$$};
\node (W1) at (3.5,0.5){$$}; \node (W2) at (1.5,1.5){$$}; \node (W3) at (3.5,2.5){$$};  \node (W4) at (1.5,3.5){$$};
\draw[thick]  (0,0) rectangle (4,4);
%edge
\draw[line width=0.05cm]  (B1)--(0,1.5); \draw[line width=0.05cm]  (W3)--(4,1.5); 
\draw[line width=0.05cm]  (B1)--(W2); \draw[line width=0.05cm]  (B1)--(0,0.5); \draw[line width=0.05cm]  (B1)--(1,0); 
\draw[line width=0.05cm]  (W2)--(B2); \draw[line width=0.05cm]  (W2)--(B3); \draw[line width=0.05cm]  (W3)--(B2);
\draw[line width=0.05cm]  (B2)--(W1); \draw[line width=0.05cm]  (B2)--(1.75,0); 
\draw[line width=0.05cm]  (W1)--(3,0); \draw[line width=0.05cm]  (W1)--(4,0.5); \draw[line width=0.05cm]  (B3)--(W4); 
\draw[line width=0.05cm]  (B3)--(0,2.5); \draw[line width=0.05cm]  (W3)--(4,2.5); \draw[line width=0.05cm]  (W3)--(B4);
\draw[line width=0.05cm]  (W4)--(B4); \draw[line width=0.05cm]  (W4)--(1,4); \draw[line width=0.05cm]  (W4)--(1.75,4); \draw[line width=0.05cm]  (B4)--(3,4);
%black
\filldraw  [ultra thick, fill=black] (0.5,0.5) circle [radius=0.16] ;\filldraw  [ultra thick, fill=black] (2.5,1.5) circle [radius=0.16] ;
\filldraw  [ultra thick, fill=black] (0.5,2.5) circle [radius=0.16] ;\filldraw  [ultra thick, fill=black] (2.5,3.5) circle [radius=0.16] ;
%white
\draw  [ultra thick,fill=white] (3.5,0.5) circle [radius=0.16] ;\draw  [ultra thick, fill=white] (1.5,1.5) circle [radius=0.16] ;
\draw  [ultra thick,fill=white] (3.5,2.5)circle [radius=0.16] ;\draw  [ultra thick, fill=white] (1.5,3.5) circle [radius=0.16] ;
\end{tikzpicture}
}}  ;

\node (Qui2) at (4,0) 
{\scalebox{0.6}{ 
\begin{tikzpicture}
%vertex
\node (B1) at (0.5,0.5){$$}; \node (B2) at (2.5,1.5){$$}; \node (B3) at (0.5,2.5){$$};  \node (B4) at (2.5,3.5){$$};
\node (W1) at (3.5,0.5){$$}; \node (W2) at (1.5,1.5){$$}; \node (W3) at (3.5,2.5){$$};  \node (W4) at (1.5,3.5){$$};
%\draw[thick]  (0,0) rectangle (4,4); 
\node (Q0a) at (0,1){$\bullet$}; \node (Q0b) at (4,1){$\bullet$}; \node (Q1a) at (0,2){$\bullet$}; \node (Q1b) at (4,2){$\bullet$}; \node (Q2a) at (1.5,0){$\bullet$};
\node (Q2b) at (1.5,4){$\bullet$}; \node (Q3a) at (2.5,0){$\bullet$}; \node (Q3b) at (2.5,4){$\bullet$}; \node (Q4) at (2,2.5){$\bullet$};
\node (Q5a) at (0,0){$\bullet$}; \node (Q5b) at (4,0){$\bullet$}; \node (Q5c) at (4,4){$\bullet$}; \node (Q5d) at (0,4){$\bullet$}; 

%\node (Q0a) at (0,1){$0$}; \node (Q0b) at (4,1){$0$}; \node (Q1a) at (0,2){$1$}; \node (Q1b) at (4,2){$1$}; \node (Q2a) at (1.5,0){$2$};
%\node (Q2b) at (1.5,4){$2$}; \node (Q3a) at (2.5,0){$3$}; \node (Q3b) at (2.5,4){$3$}; \node (Q4) at (2,2.5){$4$};
%\node (Q5a) at (0,0){$5$}; \node (Q5b) at (4,0){$5$}; \node (Q5c) at (4,4){$5$}; \node (Q5d) at (0,4){$5$};   
%edge
\draw[lightgray, line width=0.05cm]  (B1)--(0,1.5); \draw[lightgray, line width=0.05cm]  (W3)--(4,1.5); 
\draw[lightgray, line width=0.05cm]  (B1)--(W2); \draw[lightgray, line width=0.05cm]  (B1)--(0,0.5); \draw[lightgray, line width=0.05cm]  (B1)--(1,0); 
\draw[lightgray, line width=0.05cm]  (W2)--(B2); \draw[lightgray, line width=0.05cm]  (W2)--(B3); \draw[lightgray, line width=0.05cm]  (W3)--(B2);
\draw[lightgray, line width=0.05cm]  (B2)--(W1); \draw[lightgray, line width=0.05cm]  (B2)--(1.75,0); 
\draw[lightgray, line width=0.05cm]  (W1)--(3,0); \draw[lightgray, line width=0.05cm]  (W1)--(4,0.5); \draw[lightgray, line width=0.05cm]  (B3)--(W4); 
\draw[lightgray, line width=0.05cm]  (B3)--(0,2.5); \draw[lightgray, line width=0.05cm]  (W3)--(4,2.5); \draw[lightgray, line width=0.05cm]  (W3)--(B4);
\draw[lightgray, line width=0.05cm]  (W4)--(B4); \draw[lightgray, line width=0.05cm]  (W4)--(1,4); 
\draw[lightgray, line width=0.05cm]  (W4)--(1.75,4); \draw[lightgray, line width=0.05cm]  (B4)--(3,4);
%black
\filldraw  [ultra thick, draw=lightgray, fill=lightgray] (0.5,0.5) circle [radius=0.16] ;\filldraw  [ultra thick, draw=lightgray, fill=lightgray] (2.5,1.5) circle [radius=0.16] ;
\filldraw  [ultra thick, draw=lightgray, fill=lightgray] (0.5,2.5) circle [radius=0.16] ;\filldraw  [ultra thick, draw=lightgray, fill=lightgray] (2.5,3.5) circle [radius=0.16] ;
%white
\draw  [ultra thick,draw=lightgray,fill=white] (3.5,0.5) circle [radius=0.16] ;\draw  [ultra thick, draw=lightgray,fill=white] (1.5,1.5) circle [radius=0.16] ;
\draw  [ultra thick,draw=lightgray,fill=white] (3.5,2.5)circle [radius=0.16] ;\draw  [ultra thick, draw=lightgray,fill=white] (1.5,3.5) circle [radius=0.16] ;

\draw[->, line width=0.067cm] (Q0a)--node[midway,xshift=-10pt,yshift=3pt] {{\LARGE$\frac{3}{4}$}}(Q5a); 
\draw[->, line width=0.067cm] (Q0b)--node[midway,xshift=10pt,yshift=3pt] {{\LARGE$\frac{3}{4}$}}(Q5b); 
\draw[->, line width=0.067cm] (Q0b)--node[midway,xshift=5pt,yshift=8pt] {{\LARGE$\frac{1}{2}$}}(Q4); 
\draw[->, line width=0.067cm] (Q1a)--node[midway,xshift=-10pt,yshift=3pt] {{\LARGE$\frac{1}{4}$}}(Q0a); 
\draw[->, line width=0.067cm] (Q1b)--node[midway,xshift=10pt,yshift=3pt] {{\LARGE$\frac{1}{4}$}}(Q0b); 
\draw[->, line width=0.067cm] (Q1a)--node[midway,xshift=-3pt,yshift=10pt] {{\LARGE$\frac{3}{4}$}}(Q4); 
\draw[->, line width=0.067cm] (Q2a)--node[midway,xshift=-3pt,yshift=-15pt] {{\LARGE$\frac{1}{4}$}}(Q3a); 
\draw[->, line width=0.067cm] (Q2b)--node[midway,xshift=-3pt,yshift=15pt] {{\LARGE$\frac{1}{4}$}}(Q3b); 
\draw[->, line width=0.067cm] (Q2a)--node[midway,xshift=5pt,yshift=7pt] {{\LARGE$\frac{1}{2}$}}(Q1a);
\draw[->, line width=0.067cm] (Q3a)--node[midway,xshift=-3pt,yshift=10pt] {{\LARGE$\frac{1}{2}$}}(Q0b); 
\draw[->, line width=0.067cm] (Q3b)--node[midway,xshift=-7pt,yshift=3pt] {{\LARGE$\frac{1}{2}$}}(Q4);  
\draw[->, line width=0.067cm] (Q4)--node[midway,xshift=5pt,yshift=0pt] {{\LARGE$\frac{3}{4}$}}(Q2a); 
\draw[->, line width=0.067cm] (Q4)--node[midway,xshift=3pt,yshift=-10pt] {{\LARGE$\frac{3}{4}$}}(Q5c); 
\draw[->, line width=0.067cm] (Q4)--node[midway,xshift=12pt,yshift=3pt] {{\LARGE$\frac{3}{4}$}}(Q5d); 
\draw[->, line width=0.067cm] (Q5a)--node[midway,xshift=0pt,yshift=-15pt] {{\LARGE$\frac{1}{2}$}}(Q2a);   
\draw[->, line width=0.067cm] (Q5d)--node[midway,xshift=0pt,yshift=15pt] {{\LARGE$\frac{1}{2}$}}(Q2b); 
\draw[->, line width=0.067cm] (Q5b)--node[midway,xshift=0pt,yshift=-15pt] {{\LARGE$\frac{3}{4}$}}(Q3a); 
\draw[->, line width=0.067cm] (Q5c)--node[midway,xshift=0pt,yshift=15pt] {{\LARGE$\frac{3}{4}$}}(Q3b);  
\draw[->, line width=0.067cm] (Q5d)--node[midway,xshift=-10pt,yshift=0pt] {{\LARGE$\frac{1}{2}$}}(Q1a); 
\draw[->, line width=0.067cm] (Q5c)--node[midway,xshift=10pt,yshift=0pt] {{\LARGE$\frac{1}{2}$}}(Q1b); 
\end{tikzpicture}
}} ;

\end{tikzpicture}
\end{center}
Again, by the results in \cite[subsection~5.7]{Nak}, we can see rank one MCM modules arising from isoradial dimer models are actually conic, 
and non-conic divisorial ideals arise from a consistent dimer model which is not isoradial. 

%%%%%%%%%%%%%%%%%%%%%%%%%%%%%%%%%%%%%%%%%%%%%%%%%%%%%%%
\subsubsection{{\bf \textup{Type 6c}}}
\label{type6c}

We consider the reflexive polygon of type 6c. 
Thus, let $R$ be the $3$-dimensional complete local Gorenstein toric singularity defined by the cone $\sigma$: 
\[
\sigma=\mathrm{Cone}\{v_1=(1,0,1), v_2=(0,1,1), v_3=(-2,-1,1), v_4=(0,-1,1) \}. 
\]
As elements in $\Cl(R)$, we obtain $[I_1]=2[I_3]$, $[I_2]-[I_3]-[I_4]=0$, $4[I_3]+2[I_4]=0$.  
Thus, we have $\Cl(R)\cong\ZZ\times \ZZ/2\ZZ$, and each divisorial ideal is represented by $T(0,0,c,d)$ where $c\in\ZZ, d\in\ZZ/2\ZZ$. 
In this case, there are two consistent dimer models associated with $R$. 
By the results in \cite[subsection~5.8]{Nak}, we have rank one MCM $R$-modules as in Figure~\ref{6c} from such consistent dimer models. 
In this figure, each circle stands for $(c,d)\in\Cl(R)$ corresponding to an MCM module $T(0,0,c,d)$. 
Especially, a double circle stands for the origin $(0,0)$. Also, by the direct computation below, we can see black circles are conic. 
On the other hand, white circles are not conic (e.g. check the condition in Lemma~\ref{conic_check}). 

\begin{center}
\begin{tabular}{ll}
$T(0,0,0,1)\cong T(\lambda(\frac{1}{8},\frac{-3}{8},\frac{-1}{4}))$,& $T(0,0,1,0)\cong T(\lambda(\frac{-1}{4},\frac{-1}{4},\frac{-1}{2}))$,\\ 
$T(0,0,1,1)\cong T(\lambda(\frac{-1}{4},\frac{-1}{2},\frac{-1}{4}))$,& $T(0,0,2,0)\cong T(\lambda(\frac{-5}{8},\frac{-1}{8},\frac{-1}{4}))$,\\
$T(0,0,2,1)\cong T(\lambda(\frac{-1}{2},\frac{-1}{2},\frac{-1}{8}))$. 
\end{tabular} 
\end{center}

\begin{figure}[h]
{\scalebox{0.55}{
\begin{tikzpicture}
\draw [step=1,thin, gray] (-3,-1) grid (6,2);

\draw (0,0) circle [radius=0.3]; 
\filldraw (0,0) circle [radius=0.15]; 
\filldraw  [thick, fill=black] (0,1) circle [radius=0.17] ; \filldraw  [thick, fill=black] (1,0) circle [radius=0.17] ;
\filldraw  [thick, fill=black] (1,1) circle [radius=0.17] ; \filldraw  [thick, fill=black] (2,0) circle [radius=0.17] ;
\filldraw  [thick, fill=black] (2,1) circle [radius=0.17] ; \filldraw  [thick, fill=black] (3,1) circle [radius=0.17] ; 
\filldraw  [thick, fill=black] (4,1) circle [radius=0.17] ; 
\filldraw  [thick, fill=black] (-1,0) circle [radius=0.17] ; \filldraw  [thick, fill=black] (-2,0) circle [radius=0.17] ;

\filldraw  [thick, fill=white] (5,1) circle [radius=0.17] ; \filldraw  [thick, fill=white] (-1,1) circle [radius=0.17] ;
\end{tikzpicture}
}}
\caption{MCMs arising from consistent dimer models for Type 6c} 
\label{6c}
\end{figure}

As we mentioned, there are two consistent dimer models associated with $R$, 
and one of them is isoradial and the other is not isoradial. 
The following figures are an isoradial one and the associated quiver. (The fraction on each arrow is an $\sfR$-charge which indicates they are isoradial.)

\begin{center}
\begin{tikzpicture} 
\node (DM) at (0,0) 
{\scalebox{0.48}{
\begin{tikzpicture}
%vertex
\node (B1) at (1.5,0.5){$$}; \node (B2) at (4.5,2){$$}; \node (B3) at (3.5,3.5){$$};  \node (B4) at (1.5,4){$$};
\node (W1) at (4,0.5){$$}; \node (W2) at (3,2){$$}; \node (W3) at (1,2.5){$$};  \node (W4) at (0.5,4.5){$$};
\draw[thick]  (0,0) rectangle (5,5);
%edge
\draw[line width=0.06cm] (B1)--(W2); \draw[line width=0.06cm] (B1)--(W3); \draw[line width=0.06cm] (B2)--(W1); 
\draw[line width=0.06cm] (B2)--(W2); \draw[line width=0.06cm] (B3)--(W2); \draw[line width=0.06cm] (B3)--(W3); 
\draw[line width=0.06cm] (B4)--(W3); \draw[line width=0.06cm] (B4)--(W4); 
\draw[line width=0.06cm] (B1)--(0,0.5); \draw[line width=0.06cm] (W1)--(5,0.5); \draw[line width=0.06cm] (B1)--(1,0); 
\draw[line width=0.06cm] (W4)--(1,5); \draw[line width=0.06cm] (W1)--(3.16,0); \draw[line width=0.06cm] (B4)--(3.165,5); 
\draw[line width=0.06cm] (W1)--(3.875,0); \draw[line width=0.06cm] (B3)--(3.875,5); 
\draw[line width=0.06cm] (B2)--(5,2.165); \draw[line width=0.06cm] (W3)--(0,2.165); 
\draw[line width=0.06cm] (B3)--(5,4.25); \draw[line width=0.06cm] (W4)--(0,4.25);  
%black
\filldraw  [ultra thick, fill=black] (1.5,0.5) circle [radius=0.2] ;\filldraw  [ultra thick, fill=black] (4.5,2) circle [radius=0.2] ;
\filldraw  [ultra thick, fill=black] (3.5,3.5) circle [radius=0.2] ;\filldraw  [ultra thick, fill=black] (1.5,4) circle [radius=0.2] ;
%white
\draw  [ultra thick,fill=white] (4,0.5) circle [radius=0.2] ;\draw  [ultra thick, fill=white] (3,2) circle [radius=0.2] ;
\draw  [ultra thick,fill=white] (1,2.5)circle [radius=0.2] ;\draw  [ultra thick, fill=white] (0.5,4.5) circle [radius=0.2] ;
\end{tikzpicture}
}}  ;

\node (Qui) at (4,0) 
{\scalebox{0.48}{ 
\begin{tikzpicture}
%vertex
\node (B1) at (1.5,0.5){$$}; \node (B2) at (4.5,2){$$}; \node (B3) at (3.5,3.5){$$};  \node (B4) at (1.5,4){$$};
\node (W1) at (4,0.5){$$}; \node (W2) at (3,2){$$}; \node (W3) at (1,2.5){$$};  \node (W4) at (0.5,4.5){$$};
%\draw[thick]  (0,0) rectangle (5,5);
\node (Q0a) at (0,1.25){{\Large$\bullet$}}; \node (Q0b) at (5,1.25){{\Large$\bullet$}}; 
\node (Q1a) at (0,0){{\Large$\bullet$}}; \node (Q1b) at (5,0){{\Large$\bullet$}}; \node (Q1c) at (5,5){{\Large$\bullet$}}; \node (Q1d) at (0,5){{\Large$\bullet$}}; 
\node (Q2a) at (2,0){{\Large$\bullet$}}; \node (Q2b) at (2,5){{\Large$\bullet$}}; 
\node (Q3a) at (3.5,0){{\Large$\bullet$}}; \node (Q3b) at (3.5,5){{\Large$\bullet$}}; \node (Q4) at (2.25,2.25){{\Large$\bullet$}}; 
\node (Q5a) at (0,3){{\Large$\bullet$}}; \node (Q5b) at (5,3){{\Large$\bullet$}}; 

%\node (Q0a) at (0,1.25){{\Large$0$}}; \node (Q0b) at (5,1.25){{\Large$0$}}; 
%\node (Q1a) at (0,0){{\Large$1$}}; \node (Q1b) at (5,0){{\Large$1$}}; \node (Q1c) at (5,5){{\Large$1$}}; \node (Q1d) at (0,5){{\Large$1$}}; 
%\node (Q2a) at (2,0){{\Large$2$}}; \node (Q2b) at (2,5){{\Large$2$}}; 
%\node (Q3a) at (3.5,0){{\Large$3$}}; \node (Q3b) at (3.5,5){{\Large$3$}}; \node (Q4) at (2.25,2.25){{\Large$4$}}; 
%\node (Q5a) at (0,3){{\Large$5$}}; \node (Q5b) at (5,3){{\Large$5$}}; 
%edge
\draw[lightgray, line width=0.06cm] (B1)--(W2); \draw[lightgray, line width=0.06cm] (B1)--(W3); \draw[lightgray, line width=0.06cm] (B2)--(W1); 
\draw[lightgray, line width=0.06cm] (B2)--(W2); \draw[lightgray, line width=0.06cm] (B3)--(W2); \draw[lightgray, line width=0.06cm] (B3)--(W3); 
\draw[lightgray, line width=0.06cm] (B4)--(W3); \draw[lightgray, line width=0.06cm] (B4)--(W4); 
\draw[lightgray, line width=0.06cm] (B1)--(0,0.5); \draw[lightgray, line width=0.06cm] (W1)--(5,0.5); \draw[lightgray, line width=0.06cm] (B1)--(1,0); 
\draw[lightgray, line width=0.06cm] (W4)--(1,5); \draw[lightgray, line width=0.06cm] (W1)--(3.16,0); \draw[lightgray, line width=0.06cm] (B4)--(3.165,5); 
\draw[lightgray, line width=0.06cm] (W1)--(3.875,0); \draw[lightgray, line width=0.06cm] (B3)--(3.875,5); 
\draw[lightgray, line width=0.06cm] (B2)--(5,2.165); \draw[lightgray, line width=0.06cm] (W3)--(0,2.165); 
\draw[lightgray, line width=0.06cm] (B3)--(5,4.25); \draw[lightgray, line width=0.06cm] (W4)--(0,4.25);  
%black
\filldraw  [ultra thick, draw=lightgray, fill=lightgray] (1.5,0.5) circle [radius=0.2] ;\filldraw  [ultra thick, draw=lightgray, fill=lightgray] (4.5,2) circle [radius=0.2] ;
\filldraw  [ultra thick, draw=lightgray, fill=lightgray] (3.5,3.5) circle [radius=0.2] ;\filldraw  [ultra thick, draw=lightgray, fill=lightgray] (1.5,4) circle [radius=0.2] ;
%white
\draw  [ultra thick,draw=lightgray,fill=white] (4,0.5) circle [radius=0.2] ;\draw  [ultra thick,draw=lightgray,fill=white] (3,2) circle [radius=0.2] ;
\draw  [ultra thick,draw=lightgray,fill=white] (1,2.5)circle [radius=0.2] ;\draw  [ultra thick, draw=lightgray,fill=white] (0.5,4.5) circle [radius=0.2] ;
%arrow
\draw[->, line width=0.08cm] (Q0a)--node[midway,xshift=-10pt,yshift=3pt] {{\huge$\frac{1}{2}$}}(Q1a); 
\draw[->, line width=0.08cm] (Q0a)--node[midway,xshift=-10pt,yshift=0pt] {{\huge$\frac{3}{4}$}}(Q5a); 
\draw[->, line width=0.08cm] (Q0b)--node[midway,xshift=10pt,yshift=3pt] {{\huge$\frac{1}{2}$}}(Q1b); 
\draw[->, line width=0.08cm] (Q0b)--node[midway,xshift=10pt,yshift=0pt] {{\huge$\frac{3}{4}$}}(Q5b); 
\draw[->, line width=0.08cm] (Q1a)--node[midway,xshift=0pt,yshift=-15pt] {{\huge$\frac{1}{2}$}}(Q2a); 
\draw[->, line width=0.08cm] (Q1b)--node[midway,xshift=0pt,yshift=-15pt] {{\huge$\frac{1}{4}$}}(Q3a); 
\draw[->, line width=0.08cm] (Q1c)--node[midway,xshift=0pt,yshift=15pt] {{\huge$\frac{1}{4}$}}(Q3b); 
\draw[->, line width=0.08cm] (Q1d)--node[midway,xshift=0pt,yshift=15pt] {{\huge$\frac{1}{2}$}}(Q2b); 
\draw[->, line width=0.08cm] (Q2a)--node[midway,xshift=15pt,yshift=20pt] {{\huge$\frac{1}{2}$}}(Q0b); 
\draw[->, line width=0.08cm] (Q2a)--node[midway,xshift=-7pt,yshift=0pt] {{\huge$\frac{3}{4}$}}(Q4); 
\draw[->, line width=0.08cm] (Q2b)--node[midway,xshift=-5pt,yshift=10pt] {{\huge$\frac{3}{4}$}}(Q5a); 
\draw[->, line width=0.08cm] (Q3a)--node[midway,xshift=0pt,yshift=-15pt] {{\huge$\frac{3}{4}$}}(Q2a); 
\draw[->, line width=0.08cm] (Q3b)--node[midway,xshift=0pt,yshift=15pt] {{\huge$\frac{3}{4}$}}(Q2b); 
\draw[->, line width=0.08cm] (Q3b)--node[midway,xshift=8pt,yshift=0pt] {{\huge$\frac{1}{2}$}}(Q4); 
\draw[->, line width=0.08cm] (Q4)--node[midway,xshift=-5pt,yshift=12pt] {{\huge$\frac{1}{4}$}}(Q0a); 
\draw[->, line width=0.08cm] (Q4)--node[midway,xshift=8pt,yshift=17pt] {{\huge$\frac{1}{2}$}}(Q5b); 
\draw[->, line width=0.08cm] (Q5a)--node[midway,xshift=-10pt,yshift=0pt] {{\huge$\frac{3}{4}$}}(Q1d); 
\draw[->, line width=0.08cm] (Q5a)--node[midway,xshift=8pt,yshift=-10pt] {{\huge$\frac{1}{2}$}}(Q3b); 
\draw[->, line width=0.08cm] (Q5b)--node[midway,xshift=10pt,yshift=0pt] {{\huge$\frac{3}{4}$}}(Q1c); 
\draw[->, line width=0.08cm] (Q5b)--node[midway,xshift=-10pt,yshift=5pt] {{\huge$\frac{3}{4}$}}(Q2a); 
\end{tikzpicture}
}} ;

\end{tikzpicture}
\end{center}
Again, by the results in \cite[subsection~5.8]{Nak}, we can see rank one MCM modules arising from an isoradial dimer model is actually conic, 
and non-conic divisorial ideals arise from a consistent dimer model which is not isoradial. 

%%%%%%%%%%%%%%%%%%%%%%%%%%%%%%%%%%%%%%%%%%%%%%%%%%%%%%%
\subsubsection{{\bf \textup{Type 7a}}}
\label{type7a}

We consider the reflexive polygon of type 7a. 
Thus, let $R$ be the $3$-dimensional complete local Gorenstein toric singularity defined by the cone $\sigma$: 
\[
\sigma=\mathrm{Cone}\{v_1=(1,0,1), v_2=(0,1,1), v_3=(-1,1,1), v_4=(-1,-1,1), v_5=(1,-1,1) \}. 
\]
As elements in $\Cl(R)$, we obtain $[I_1]+2[I_4]+2[I_5]=0$, $[I_2]-4[I_4]-2[I_5]=0$, $[I_3]+3[I_4]+[I_5]=0$. 
Therefore, we have $\Cl(R)\cong\ZZ^2$, and each divisorial ideal is represented by $T(0,0,0,d,e)$ where $d, e\in\ZZ$. 
In this case, there are four consistent dimer models associated with $R$. 
By the results in \cite[subsection~5.9]{Nak}, we have rank one MCM $R$-modules as in Figure~\ref{7a} from such consistent dimer models. 
In the following figure, each circle stands for $(d,e)\in\Cl(R)$ corresponding to an MCM module $T(0,0,0,d,e)$. 
Especially, a double circle stands for the origin $(0,0)$. Also, by the direct computation below, we can see black circles are conic. 
On the other hand, white circles are not conic (e.g. check the condition in Lemma~\ref{conic_check}). 

\begin{center}
\begin{tabular}{ll}
$T(0,0,0,0,1)\cong T(\lambda(\frac{1}{8},\frac{-1}{8},\frac{-1}{8}))$,& $T(0,0,0,1,0)\cong T(\lambda(\frac{-1}{2},\frac{-1}{2},\frac{-3}{8}))$,\\ 
$T(0,0,0,1,1)\cong T(\lambda(\frac{-1}{4},\frac{-3}{4},\frac{-1}{8}))$,& $T(0,0,0,2,0)\cong T(\lambda(\frac{-3}{4},\frac{-3}{4},\frac{-1}{8}))$, \\
$T(0,0,0,2,1)\cong T(\lambda(\frac{-3}{8},\frac{-7}{8},\frac{1}{4}))$,& $T(0,0,0,2,2)\cong T(\lambda(\frac{-1}{4},\frac{-9}{8},\frac{1}{4}))$, \\
$T(0,0,0,3,1)\cong T(\lambda(\frac{-3}{4},\frac{-5}{4},\frac{3}{8}))$,& $T(0,0,0,3,2)\cong T(\lambda(\frac{-5}{8},\frac{-11}{8},\frac{1}{2}))$,\\ 
$T(0,0,0,4,2)\cong T(\lambda(\frac{-7}{8},\frac{-25}{16},\frac{5}{8}))$. 
\end{tabular} 
\end{center}

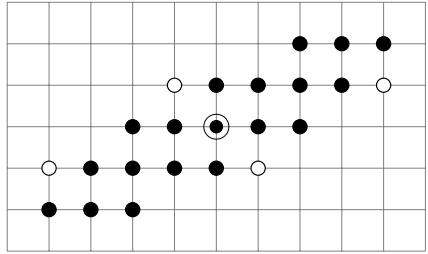
\begin{figure}[h]
{\scalebox{0.55}{
\begin{tikzpicture}
\draw [step=1,thin, gray] (-5,-3) grid (5,3);

\draw (0,0) circle [radius=0.3]; 
\filldraw (0,0) circle [radius=0.15]; 
\filldraw  [thick, fill=black] (0,1) circle [radius=0.17] ; \filldraw  [thick, fill=black] (1,0) circle [radius=0.17] ;
\filldraw  [thick, fill=black] (1,1) circle [radius=0.17] ; \filldraw  [thick, fill=black] (2,0) circle [radius=0.17] ;
\filldraw  [thick, fill=black] (2,1) circle [radius=0.17] ; \filldraw  [thick, fill=black] (3,1) circle [radius=0.17] ; 
\filldraw  [thick, fill=black] (2,2) circle [radius=0.17] ; \filldraw  [thick, fill=black] (3,2) circle [radius=0.17] ;
\filldraw  [thick, fill=black] (4,2) circle [radius=0.17] ; 
\filldraw  [thick, fill=black] (0,-1) circle [radius=0.17] ; \filldraw  [thick, fill=black] (-1,0) circle [radius=0.17] ;
\filldraw  [thick, fill=black] (-1,-1) circle [radius=0.17] ; \filldraw  [thick, fill=black] (-2,0) circle [radius=0.17] ;
\filldraw  [thick, fill=black] (-2,-1) circle [radius=0.17] ; \filldraw  [thick, fill=black] (-3,-1) circle [radius=0.17] ; 
\filldraw  [thick, fill=black] (-2,-2) circle [radius=0.17] ; \filldraw  [thick, fill=black] (-3,-2) circle [radius=0.17] ;
\filldraw  [thick, fill=black] (-4,-2) circle [radius=0.17] ; 

\filldraw  [thick, fill=white] (1,-1) circle [radius=0.17] ; \filldraw  [thick, fill=white] (-1,1) circle [radius=0.17] ;
\filldraw  [thick, fill=white] (4,1) circle [radius=0.17] ; \filldraw  [thick, fill=white] (-4,-1) circle [radius=0.17] ;
\end{tikzpicture}
}}
\caption{MCMs arising from consistent dimer models for Type 7a} 
\label{7a}
\end{figure}

As we mentioned, there are four consistent dimer models associated with $R$, 
and one of them is isoradial and the others are not isoradial. 
The following figures are an isoradial one and the associated quiver. (The fraction on each arrow is an $\sfR$-charge which indicates this dimer model is isoradial.)

%%%%%%%%%%%%%%%%%%%%%%%%%%%%%%%%%%%
%%%%%%%%%%%%Type_7a-1%%%%%%%%%%%%%%%%%%
%%%%%%%%%%%%%%%%%%%%%%%%%%%%%%%%%%%

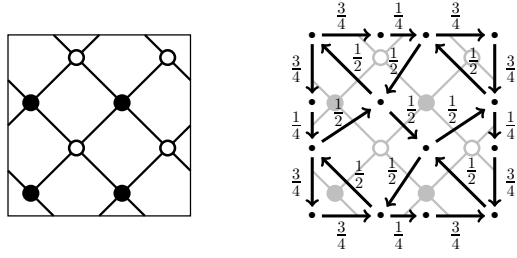
\begin{center}
\begin{tikzpicture} 
\node (DM) at (0,0) 
{\scalebox{0.6}{
\begin{tikzpicture}
%vertex
\node (B1) at (0.5,0.5){$$}; \node (B2) at (2.5,0.5){$$}; \node (B3) at (0.5,2.5){$$};  \node (B4) at (2.5,2.5){$$};
\node (W1) at (1.5,1.5){$$}; \node (W2) at (3.5,1.5){$$}; \node (W3) at (1.5,3.5){$$};  \node (W4) at (3.5,3.5){$$};
\draw[thick]  (0,0) rectangle (4,4);

%edge
\draw[line width=0.05cm]  (B1)--(W1); \draw[line width=0.05cm]  (B1)--(1,0); \draw[line width=0.05cm]  (B1)--(0,1);
\draw[line width=0.05cm]  (B2)--(W1); \draw[line width=0.05cm]  (B2)--(W2); \draw[line width=0.05cm]  (B2)--(2,0); 
\draw[line width=0.05cm]  (B2)--(3,0); \draw[line width=0.05cm]  (B3)--(W1); \draw[line width=0.05cm]  (B3)--(W3); 
\draw[line width=0.05cm]  (B3)--(0,2); \draw[line width=0.05cm]  (B3)--(0,3); \draw[line width=0.05cm]  (B4)--(W1);
\draw[line width=0.05cm]  (B4)--(W2); \draw[line width=0.05cm]  (B4)--(W3); \draw[line width=0.05cm]  (B4)--(W4); 
\draw[line width=0.05cm]  (W2)--(4,1); \draw[line width=0.05cm]  (W2)--(4,2); \draw[line width=0.05cm]  (W3)--(1,4); 
\draw[line width=0.05cm]  (W3)--(2,4); \draw[line width=0.05cm]  (W4)--(4,3); \draw[line width=0.05cm]  (W4)--(3,4);
%black
\filldraw  [ultra thick, fill=black] (0.5,0.5) circle [radius=0.16] ;\filldraw  [ultra thick, fill=black] (2.5,0.5) circle [radius=0.16] ;
\filldraw  [ultra thick, fill=black] (0.5,2.5) circle [radius=0.16] ;\filldraw  [ultra thick, fill=black] (2.5,2.5) circle [radius=0.16] ;
%white
\draw  [ultra thick,fill=white] (1.5,1.5) circle [radius=0.16] ;\draw  [ultra thick, fill=white] (3.5,1.5) circle [radius=0.16] ;
\draw  [ultra thick,fill=white] (1.5,3.5)circle [radius=0.16] ;\draw  [ultra thick, fill=white] (3.5,3.5) circle [radius=0.16] ;
\end{tikzpicture}
}}  ;

\node (Qui) at (4,0) 
{\scalebox{0.6}{ 
\begin{tikzpicture}
%vertex
\node (B1) at (0.5,0.5){$$}; \node (B2) at (2.5,0.5){$$}; \node (B3) at (0.5,2.5){$$};  \node (B4) at (2.5,2.5){$$};
\node (W1) at (1.5,1.5){$$}; \node (W2) at (3.5,1.5){$$}; \node (W3) at (1.5,3.5){$$};  \node (W4) at (3.5,3.5){$$};
%\draw[thick]  (0,0) rectangle (4,4);
\node (Q0a) at (0,2.5){$\bullet$}; \node (Q0b) at (4,2.5){$\bullet$}; \node (Q1a) at (0,1.5){$\bullet$}; \node (Q1b) at (4,1.5){$\bullet$};
\node (Q2a) at (0,0){$\bullet$}; \node (Q2b) at (4,0){$\bullet$}; \node (Q2c) at (4,4){$\bullet$}; \node (Q2d) at (0,4){$\bullet$}; 
\node (Q3) at (2.5,1.5){$\bullet$}; \node (Q4a) at (1.5,0){$\bullet$}; \node (Q4b) at (1.5,4){$\bullet$};  
\node (Q5) at (1.5,2.5){$\bullet$}; \node (Q6a) at (2.5,0){$\bullet$}; \node (Q6b) at (2.5,4){$\bullet$};  
   
%\node (Q0a) at (0,2.5){$0$}; \node (Q0b) at (4,2.5){$0$}; \node (Q1a) at (0,1.5){$1$}; \node (Q1b) at (4,1.5){$1$};
%\node (Q2a) at (0,0){$2$}; \node (Q2b) at (4,0){$2$}; \node (Q2c) at (4,4){$2$}; \node (Q2d) at (0,4){$2$}; 
%\node (Q3) at (2.5,1.5){$3$}; \node (Q4a) at (1.5,0){$4$}; \node (Q4b) at (1.5,4){$4$};  
%\node (Q5) at (1.5,2.5){$5$}; \node (Q6a) at (2.5,0){$6$}; \node (Q6b) at (2.5,4){$6$}; 
%edge
\draw[lightgray, line width=0.05cm]  (B1)--(W1); \draw[lightgray, line width=0.05cm]  (B1)--(1,0); \draw[lightgray, line width=0.05cm]  (B1)--(0,1);
\draw[lightgray, line width=0.05cm]  (B2)--(W1); \draw[lightgray, line width=0.05cm]  (B2)--(W2); \draw[lightgray, line width=0.05cm]  (B2)--(2,0); 
\draw[lightgray, line width=0.05cm]  (B2)--(3,0); \draw[lightgray, line width=0.05cm]  (B3)--(W1); \draw[lightgray, line width=0.05cm]  (B3)--(W3); 
\draw[lightgray, line width=0.05cm]  (B3)--(0,2); \draw[lightgray, line width=0.05cm]  (B3)--(0,3); \draw[lightgray, line width=0.05cm]  (B4)--(W1);
\draw[lightgray, line width=0.05cm]  (B4)--(W2); \draw[lightgray, line width=0.05cm]  (B4)--(W3); \draw[lightgray, line width=0.05cm]  (B4)--(W4); 
\draw[lightgray, line width=0.05cm]  (W2)--(4,1); \draw[lightgray, line width=0.05cm]  (W2)--(4,2); \draw[lightgray, line width=0.05cm]  (W3)--(1,4); 
\draw[lightgray, line width=0.05cm]  (W3)--(2,4); \draw[lightgray, line width=0.05cm]  (W4)--(4,3); \draw[lightgray, line width=0.05cm]  (W4)--(3,4);
%black
\filldraw  [ultra thick,draw=lightgray, fill=lightgray] (0.5,0.5) circle [radius=0.16] ;\filldraw  [ultra thick, draw=lightgray, fill=lightgray] (2.5,0.5) circle [radius=0.16] ;
\filldraw  [ultra thick, draw=lightgray, fill=lightgray] (0.5,2.5) circle [radius=0.16] ;\filldraw  [ultra thick, draw=lightgray, fill=lightgray] (2.5,2.5) circle [radius=0.16] ;
%white
\draw  [ultra thick,draw=lightgray,fill=white] (1.5,1.5) circle [radius=0.16] ;\draw  [ultra thick, draw=lightgray,fill=white] (3.5,1.5) circle [radius=0.16] ;
\draw  [ultra thick,draw=lightgray,fill=white] (1.5,3.5)circle [radius=0.16] ;\draw  [ultra thick, draw=lightgray,fill=white] (3.5,3.5) circle [radius=0.16] ;
%arrow
\draw[->, line width=0.067cm] (Q0a)--node[midway,xshift=-10pt,yshift=0pt] {{\LARGE$\frac{1}{4}$}}(Q1a); 
\draw[->, line width=0.067cm] (Q0b)--node[midway,xshift=10pt,yshift=0pt] {{\LARGE$\frac{1}{4}$}}(Q1b);   
\draw[->, line width=0.067cm] (Q1a)--node[midway,xshift=-10pt,yshift=0pt] {{\LARGE$\frac{3}{4}$}}(Q2a); 
\draw[->, line width=0.067cm] (Q1b)--node[midway,xshift=10pt,yshift=0pt] {{\LARGE$\frac{3}{4}$}}(Q2b); 
\draw[->, line width=0.067cm] (Q2a)--node[midway,xshift=-3pt,yshift=-12pt] {{\LARGE$\frac{3}{4}$}}(Q4a); 
\draw[->, line width=0.067cm] (Q2c)--node[midway,xshift=10pt,yshift=0pt] {{\LARGE$\frac{3}{4}$}}(Q0b); 
\draw[->, line width=0.067cm] (Q2d)--node[midway,xshift=-10pt,yshift=0pt] {{\LARGE$\frac{3}{4}$}}(Q0a); 
\draw[->, line width=0.067cm] (Q2d)--node[midway,xshift=-3pt,yshift=12pt] {{\LARGE$\frac{3}{4}$}}(Q4b); 
\draw[->, line width=0.067cm] (Q4a)--node[midway,xshift=-3pt,yshift=-12pt] {{\LARGE$\frac{1}{4}$}}(Q6a); 
\draw[->, line width=0.067cm] (Q4b)--node[midway,xshift=-3pt,yshift=12pt] {{\LARGE$\frac{1}{4}$}}(Q6b); 
\draw[->, line width=0.067cm] (Q6a)--node[midway,xshift=-3pt,yshift=-12pt] {{\LARGE$\frac{3}{4}$}}(Q2b); 
\draw[->, line width=0.067cm] (Q6b)--node[midway,xshift=-3pt,yshift=12pt] {{\LARGE$\frac{3}{4}$}}(Q2c); 
\draw[->, line width=0.067cm] (Q0b)--node[midway,xshift=8pt,yshift=5pt] {{\LARGE$\frac{1}{2}$}}(Q6b);  
\draw[->, line width=0.067cm] (Q1a)--node[midway,xshift=-5pt,yshift=8pt] {{\LARGE$\frac{1}{2}$}}(Q5); 
\draw[->, line width=0.067cm] (Q2b)--node[midway,xshift=7pt,yshift=7pt] {{\LARGE$\frac{1}{2}$}}(Q3); 
\draw[->, line width=0.067cm] (Q3)--node[midway,xshift=-7pt,yshift=7pt] {{\LARGE$\frac{1}{2}$}}(Q4a); 
\draw[->, line width=0.067cm] (Q3)--node[midway,xshift=-5pt,yshift=10pt] {{\LARGE$\frac{1}{2}$}}(Q0b); 
\draw[->, line width=0.067cm] (Q4a)--node[midway,xshift=8pt,yshift=5pt] {{\LARGE$\frac{1}{2}$}}(Q1a);
\draw[->, line width=0.067cm] (Q5)--node[midway,xshift=7pt,yshift=7pt] {{\LARGE$\frac{1}{2}$}}(Q2d);
\draw[->, line width=0.067cm] (Q5)--node[midway,xshift=5pt,yshift=10pt] {{\LARGE$\frac{1}{2}$}}(Q3); 
\draw[->, line width=0.067cm] (Q6b)--node[midway,xshift=-6pt,yshift=5pt] {{\LARGE$\frac{1}{2}$}}(Q5);
\end{tikzpicture}
}} ;

\end{tikzpicture}
\end{center}
Again, by the results in \cite[subsection~5.9]{Nak}, we can see rank one MCM modules arising from an isoradial dimer model is actually conic, 
and non-conic divisorial ideals arise from consistent dimer models which are not isoradial. 

%%%%%%%%%%%%%%%%%%%%%%%%%%%%%%%%%%%%%%%%%%%%%%%%%%%%%%%
\subsubsection{{\bf \textup{Type 7b}}}
\label{type7b}

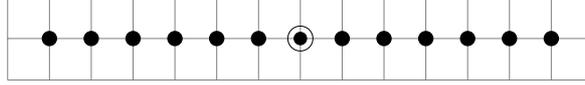
\begin{figure}[h]
{\scalebox{0.55}{
\begin{tikzpicture}
\draw [step=1,thin, gray] (-7,-1) grid (7,1);

\draw (0,0) circle [radius=0.3]; 
\filldraw (0,0) circle [radius=0.15]; 
\filldraw  [thick, fill=black] (1,0) circle [radius=0.17] ; \filldraw  [thick, fill=black] (2,0) circle [radius=0.17] ;
\filldraw  [thick, fill=black] (3,0) circle [radius=0.17] ; \filldraw  [thick, fill=black] (4,0) circle [radius=0.17] ;
\filldraw  [thick, fill=black] (5,0) circle [radius=0.17] ; \filldraw  [thick, fill=black] (6,0) circle [radius=0.17] ;
\filldraw  [thick, fill=black] (-1,0) circle [radius=0.17] ; \filldraw  [thick, fill=black] (-2,0) circle [radius=0.17] ;
\filldraw  [thick, fill=black] (-3,0) circle [radius=0.17] ; \filldraw  [thick, fill=black] (-4,0) circle [radius=0.17] ;
\filldraw  [thick, fill=black] (-5,0) circle [radius=0.17] ; \filldraw  [thick, fill=black] (-6,0) circle [radius=0.17] ;
\end{tikzpicture}
}}
\caption{MCMs arising from consistent dimer models for Type 7b} 
\label{7b}
\end{figure}

We consider the reflexive polygon of type 7b. 
Thus, let $R$ be the $3$-dimensional complete local Gorenstein toric singularity defined by the cone $\sigma$: 
\[
\sigma=\mathrm{Cone}\{v_1=(1,0,1), v_2=(0,1,1), v_3=(-2,-1,1), v_4=(1,-1,1) \}. 
\]
As elements in $\Cl(R)$, we obtain $[I_1]-6[I_3]=0$, $[I_2]+3[I_3]=0$, $[I_4]+4[I_3]=0$.
Thus, we have $\Cl(R)\cong\ZZ$, and each divisorial ideal is represented by $T(0,0,c,0)$ where $c\in\ZZ$. 
There is a unique consistent dimer model associated with $R$. 
Since the existence of isoradial dimer model for each $3$-dimensional Gorenstein toric singularity is guaranteed by \cite{Gul}, 
such a unique consistent dimer model is isoradial. 
By the results in \cite[subsection~5.10]{Nak}, we have rank one MCM $R$-modules as in Figure~\ref{7b} from such a consistent dimer model. 
In this figure, each circle stands for $c\in\Cl(R)$ corresponding to an MCM module $T(0,0,c,0)$. 
Especially, a double circle stands for the origin $(0,0)$. 
Also, by the direct computation below, we can see all rank one MCM modules arising from an isoradial dimer model are conic. 

\begin{center}
\begin{tabular}{ll}
$T(0,0,1,0)\cong T(\lambda(\frac{-1}{4},\frac{-1}{4},\frac{-1}{2}))$,& $T(0,0,2,0)\cong T(\lambda(\frac{-1}{2},\frac{-1}{2},\frac{-3}{8}))$, \\ 
$T(0,0,3,0)\cong T(\lambda(\frac{-7}{8},\frac{-1}{2},\frac{1}{4}))$, & $T(0,0,4,0)\cong T(\lambda(\frac{-5}{4},\frac{-1}{2},\frac{3}{8}))$,\\
$T(0,0,5,0)\cong T(\lambda(\frac{-13}{8},\frac{-7}{8},\frac{3}{4}))$,& $T(0,0,6,0)\cong T(\lambda(\frac{-29}{16},\frac{-7}{8},\frac{7}{8}))$. 
\end{tabular}
\end{center}

\subsubsection{{\bf \textup{Type 8a}}}
\label{type8a}

We consider the reflexive polygon of type 8a. 
Thus, let $R$ be the $3$-dimensional complete local Gorenstein toric singularity defined by the cone $\sigma$: 
\[
\sigma=\mathrm{Cone}\{v_1=(1,1,1), v_2=(-1,1,1), v_3=(-1,-1,1), v_4=(1,-1,1) \}. 
\]
As elements in $\Cl(R)$, we obtain $2[I_1]=-2[I_2]=2[I_3]=-2[I_4]$, and $[I_4]=[I_1]+[I_2]-[I_3]$.  
Therefore, we have $\Cl(R)\cong\ZZ\times(\ZZ/2\ZZ)^2$, and each divisorial ideal is represented by $T(a,b,c,0)$ where $a\in\ZZ, \, b,c\in\ZZ/2\ZZ$. 
In this case, there are four consistent dimer models associated with $R$. 
By the results in \cite[subsection~5.11]{Nak}, rank one MCM $R$-modules arising from such consistent dimer models are 
modules corresponding to elements $(0,0,0,0)$, $(1,0,0,0)$, $(0,1,0,0)$, $(0,0,1,0)$, $(1,1,0,0)$, 
$(0,1,1,0)$, $(1,1,1,0)$, $(-1,0,1,0)$, $(2,0,0,0)$, $(2,1,1,0)$, $(3,1,0,0)$ in $\Cl(R)$ and their $R$-duals. 
Also, we can see MCM modules listed below and their $R$-duals are non-free conic divisorial ideals.

\medskip

\begin{center}
\begin{tabular}{ll}
$T(1,0,0,0)\cong T(\lambda(\frac{1}{4},\frac{1}{4},\frac{-1}{4}))$,& $T(0,1,0,0)\cong T(\lambda(\frac{-1}{4},\frac{1}{4},\frac{-1}{4}))$, \\
$T(0,0,1,0)\cong T(\lambda(\frac{-1}{4},\frac{-1}{4},\frac{-1}{4}))$,& $T(1,1,0,0)\cong T(\lambda(0,\frac{1}{2},\frac{-1}{4}))$, \\
$T(0,1,1,0)\cong T(\lambda(\frac{-1}{2},0,\frac{-1}{4}))$,& $T(1,1,1,0)\cong T(\lambda(\frac{-1}{4},\frac{1}{4},\frac{1}{4}))$, \\
$T(-1,0,1,0)\cong T(\lambda(\frac{-1}{2},\frac{-1}{2},\frac{-1}{4}))$. 
\end{tabular} 
\end{center}

On the other hand, we can check remaining ones are not conic (e.g. check the condition in Lemma~\ref{conic_check}). 

As we mentioned, there are four consistent dimer models associated with $R$, 
and one of them is isoradial and the others are not isoradial. 
The following figures are an isoradial one and the associated quiver. (The fraction on each arrow is an $\sfR$-charge which indicates this dimer model is isoradial.)

%%%%%%%%%%%%%%%%%%%%%%%%%%%%%%%%%%%
%%%%%%%%%%%%Type_8a-1%%%%%%%%%%%%%%%%%%
%%%%%%%%%%%%%%%%%%%%%%%%%%%%%%%%%%%

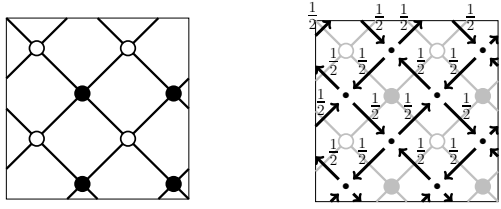
\begin{center}
\begin{tikzpicture} 
\node (DM) at (0,0) 
{\scalebox{0.4}{
\begin{tikzpicture}
%vertex
\node (B1) at (2.5,0.5){$$}; \node (B2) at (5.5,0.5){$$}; \node (B3) at (2.5,3.5){$$};  \node (B4) at (5.5,3.5){$$};
\node (W1) at (1,2){$$}; \node (W2) at (4,2){$$}; \node (W3) at (1,5){$$};  \node (W4) at (4,5){$$};
\draw[thick]  (0,0) rectangle (6,6);
%edge
\draw[line width=0.075cm] (B1)--(W1); \draw[line width=0.075cm] (B1)--(W2); \draw[line width=0.075cm] (B2)--(W2); 
\draw[line width=0.075cm] (B3)--(W1); \draw[line width=0.075cm] (B3)--(W2); \draw[line width=0.075cm] (B3)--(W3);
\draw[line width=0.075cm] (B3)--(W4); \draw[line width=0.075cm] (B4)--(W2); \draw[line width=0.075cm] (B4)--(W4); 
\draw[line width=0.075cm] (B1)--(2,0); \draw[line width=0.075cm] (B1)--(3,0); \draw[line width=0.075cm] (W3)--(2,6); \draw[line width=0.075cm] (W4)--(3,6); 
\draw[line width=0.075cm] (B2)--(5,0); \draw[line width=0.075cm] (B2)--(6,0); \draw[line width=0.075cm] (W4)--(5,6); \draw[line width=0.075cm] (W3)--(0,6); 
\draw[line width=0.075cm] (W1)--(0,1); \draw[line width=0.075cm] (W1)--(0,3); \draw[line width=0.075cm] (B2)--(6,1); \draw[line width=0.075cm] (B4)--(6,3); 
\draw[line width=0.075cm] (W3)--(0,4); \draw[line width=0.075cm] (B4)--(6,4); 
%black
\filldraw  [ultra thick, fill=black] (2.5,0.5) circle [radius=0.24] ;\filldraw  [ultra thick, fill=black] (5.5,0.5) circle [radius=0.24] ;
\filldraw  [ultra thick, fill=black] (2.5,3.5) circle [radius=0.24] ;\filldraw  [ultra thick, fill=black] (5.5,3.5) circle [radius=0.24] ;
%white
\draw  [ultra thick,fill=white] (1,2) circle [radius=0.24] ;\draw  [ultra thick, fill=white] (4,2) circle [radius=0.24] ;
\draw  [ultra thick,fill=white] (1,5)circle [radius=0.24] ;\draw  [ultra thick, fill=white] (4,5) circle [radius=0.24] ;
\end{tikzpicture}
}}  ;

\node (Qui) at (4,0.1) 
{\scalebox{0.4}{ 
\begin{tikzpicture}
%vertex
\node (B1) at (2.5,0.5){$$}; \node (B2) at (5.5,0.5){$$}; \node (B3) at (2.5,3.5){$$};  \node (B4) at (5.5,3.5){$$};
\node (W1) at (1,2){$$}; \node (W2) at (4,2){$$}; \node (W3) at (1,5){$$};  \node (W4) at (4,5){$$};
\node (Q0) at (4,3.5){{\LARGE$\bullet$}}; \node (Q1) at (2.5,2){{\LARGE$\bullet$}}; \node (Q2) at (1,3.5){{\LARGE$\bullet$}}; 
\node (Q3) at (2.5,5){{\LARGE$\bullet$}}; \node (Q4) at (4,0.5){{\LARGE$\bullet$}}; \node (Q5) at (5.5,5){{\LARGE$\bullet$}}; 
\node (Q6) at (1,0.5){{\LARGE$\bullet$}}; \node (Q7) at (5.5,2){{\LARGE$\bullet$}};

%\node (Q0) at (4,3.5){{\LARGE$0$}}; \node (Q1) at (2.5,2){{\LARGE$1$}}; \node (Q2) at (1,3.5){{\LARGE$2$}}; 
%\node (Q3) at (2.5,5){{\LARGE$3$}}; \node (Q4) at (4,0.5){{\LARGE$4$}}; \node (Q5) at (5.5,5){{\LARGE$5$}}; 
%\node (Q6) at (1,0.5){{\LARGE$6$}}; \node (Q7) at (5.5,2){{\LARGE$7$}};
\draw[thick]  (0,0) rectangle (6,6);

%edge
\draw[lightgray, line width=0.075cm] (B1)--(W1); \draw[lightgray, line width=0.075cm] (B1)--(W2); \draw[lightgray, line width=0.075cm] (B2)--(W2); 
\draw[lightgray, line width=0.075cm] (B3)--(W1); \draw[lightgray, line width=0.075cm] (B3)--(W2); \draw[lightgray, line width=0.075cm] (B3)--(W3);
\draw[lightgray, line width=0.075cm] (B3)--(W4); \draw[lightgray, line width=0.075cm] (B4)--(W2); \draw[lightgray, line width=0.075cm] (B4)--(W4); 
\draw[lightgray, line width=0.075cm] (B1)--(2,0); \draw[lightgray, line width=0.075cm] (B1)--(3,0); \draw[lightgray, line width=0.075cm] (W3)--(2,6); 
\draw[lightgray, line width=0.075cm] (W4)--(3,6); \draw[lightgray, line width=0.075cm] (B2)--(5,0); \draw[lightgray, line width=0.075cm] (B2)--(6,0); 
\draw[lightgray, line width=0.075cm] (W4)--(5,6); \draw[lightgray, line width=0.075cm] (W3)--(0,6); \draw[lightgray, line width=0.075cm] (W1)--(0,1); 
\draw[lightgray, line width=0.075cm] (W1)--(0,3); \draw[lightgray, line width=0.075cm] (B2)--(6,1); \draw[lightgray, line width=0.075cm] (B4)--(6,3); 
\draw[lightgray, line width=0.075cm] (W3)--(0,4); \draw[lightgray, line width=0.075cm] (B4)--(6,4); 
%black
\filldraw  [ultra thick, draw=lightgray, fill=lightgray] (2.5,0.5) circle [radius=0.24] ;\filldraw  [ultra thick, draw=lightgray, fill=lightgray] (5.5,0.5) circle [radius=0.24] ;
\filldraw  [ultra thick, draw=lightgray, fill=lightgray] (2.5,3.5) circle [radius=0.24] ;\filldraw  [ultra thick, draw=lightgray, fill=lightgray] (5.5,3.5) circle [radius=0.24] ;
%white
\draw  [ultra thick,draw=lightgray,fill=white] (1,2) circle [radius=0.24] ;\draw  [ultra thick, draw=lightgray,fill=white] (4,2) circle [radius=0.24] ;
\draw  [ultra thick,draw=lightgray,fill=white] (1,5)circle [radius=0.24] ;\draw  [ultra thick, draw=lightgray,fill=white] (4,5) circle [radius=0.24] ;
%arrow
\draw[->, line width=0.1cm] (Q0)--node[midway,xshift=6pt,yshift=12pt] {{\Huge$\frac{1}{2}$}}(Q3); 
\draw[->, line width=0.1cm] (Q0)--node[midway,xshift=6pt,yshift=12pt] {{\Huge$\frac{1}{2}$}}(Q7); 
\draw[->, line width=0.1cm] (Q1)--node[midway,xshift=-6pt,yshift=12pt] {{\Huge$\frac{1}{2}$}}(Q0); 
\draw[->, line width=0.1cm] (Q1)--node[midway,xshift=-6pt,yshift=12pt] {{\Huge$\frac{1}{2}$}}(Q6); 
\draw[->, line width=0.1cm] (Q2)--node[midway,xshift=6pt,yshift=12pt] {{\Huge$\frac{1}{2}$}}(Q1); 
\draw[->, line width=0.1cm] (Q2)--node[midway,xshift=6pt,yshift=15pt] {{\Huge$\frac{1}{2}$}}(0,4.5); 
\draw[->, line width=0.1cm] (6,4.5)--(Q5); 
\draw[->, line width=0.1cm] (Q3)--node[midway,xshift=-6pt,yshift=12pt] {{\Huge$\frac{1}{2}$}}(Q2); 
\draw[->, line width=0.1cm] (Q3)--node[midway,xshift=-6pt,yshift=14pt] {{\Huge$\frac{1}{2}$}}(3.5,6); 
\draw[->, line width=0.1cm] (3.5,0)--(Q4); 
\draw[->, line width=0.1cm] (Q4)--node[midway,xshift=6pt,yshift=12pt] {{\Huge$\frac{1}{2}$}}(Q1); 
\draw[->, line width=0.1cm] (Q4)--(4.5,0); 
\draw[->, line width=0.1cm] (4.5,6)--node[midway,xshift=6pt,yshift=14pt] {{\Huge$\frac{1}{2}$}}(Q5); 
\draw[->, line width=0.1cm] (Q5)--node[midway,xshift=-6pt,yshift=12pt] {{\Huge$\frac{1}{2}$}}(Q0); 
\draw[->, line width=0.1cm] (Q5)--(6,5.5); 
\draw[->, line width=0.1cm] (0,5.5)--node[midway,xshift=-10pt,yshift=12pt] {{\Huge$\frac{1}{2}$}}(0.5,6); 
\draw[->, line width=0.1cm] (0.5,0)--(Q6); 
\draw[->, line width=0.1cm] (6,1.5)--(Q7); 
\draw[->, line width=0.1cm] (Q6)--node[midway,xshift=6pt,yshift=15pt] {{\Huge$\frac{1}{2}$}}(0,1.5); 
\draw[->, line width=0.1cm] (Q6)--(1.5,0); 
\draw[->, line width=0.1cm] (1.5,6)--node[midway,xshift=6pt,yshift=14pt] {{\Huge$\frac{1}{2}$}}(Q3); 
\draw[->, line width=0.1cm] (Q7)--node[midway,xshift=-6pt,yshift=12pt] {{\Huge$\frac{1}{2}$}}(Q4); 
\draw[->, line width=0.1cm] (Q7)--(6,2.5); 
\draw[->, line width=0.1cm] (0,2.5)--node[midway,xshift=-6pt,yshift=12pt] {{\Huge$\frac{1}{2}$}}(Q2);
\end{tikzpicture}
}} ;

\end{tikzpicture}
\end{center}
Again, by the results in \cite[subsection~5.11]{Nak}, we can see rank one MCM modules arising from an isoradial dimer model is actually conic, 
and non-conic divisorial ideals arise from consistent dimer models which are not isoradial. 

%%%%%%%%%%%%%%%%%%%%%%%%%%%%%%%%%%%%%%%%%%%%%%%%%%%%%%%
\subsubsection{{\bf \textup{Type 8b}}}
\label{type8b}

We consider the reflexive polygon of type 8b. 
Thus, let $R$ be the $3$-dimensional complete local Gorenstein toric singularity defined by the cone $\sigma$: 
\[
\sigma=\mathrm{Cone}\{v_1=(0,1,1), v_2=(-1,1,1), v_3=(-1,-1,1), v_4=(2,-1,1) \}. 
\]
As elements in $\Cl(R)$, we obtain $[I_1]-2[I_3]+[I_4]=0$, $[I_2]+[I_3]-2[I_4]=0$, $2[I_3]+2[I_4]=0$. 
Thus, we have $\Cl(R)\cong\ZZ\times\ZZ/2\ZZ$, and each divisorial ideal is represented by $T(0,0,c,d)$ where $c\in\ZZ, d\in\ZZ/2\ZZ$. 
In this case, there are two consistent dimer models associated with $R$. 
By the results in \cite[subsection~5.12]{Nak}, we have rank one MCM $R$-modules as in Figure~\ref{8b} from such consistent dimer models. 
In the following figure, each circle stands for $(c,d)\in\Cl(R)$ corresponding to an MCM module $T(0,0,c,d)$. 
Especially, a double circle stands for the origin $(0,0)$. Also, by the direct computation below, we can see black circles are conic. 
On the other hand, white circles are not conic (e.g. check the condition in Lemma~\ref{conic_check}). 

\begin{center}
\begin{tabular}{ll}
$T(0,0,1,0)\cong T(\lambda(\frac{-1}{2},\frac{-1}{2},\frac{-1}{4}))$,& $T(0,0,1,1)\cong T(\lambda(0,\frac{-1}{2},\frac{-1}{8}))$,\\ 
$T(0,0,2,0)\cong T(\lambda(\frac{-5}{8},\frac{-5}{8},\frac{-1}{8}))$,& $T(0,0,2,1)\cong T(\lambda(\frac{-1}{4},\frac{-3}{4},\frac{1}{4}))$, \\
$T(0,0,3,0)\cong T(\lambda(\frac{-7}{8},\frac{-9}{8},\frac{1}{4}))$,& $T(0,0,3,1)\cong T(\lambda(\frac{-5}{8},\frac{-5}{4},\frac{3}{8}))$,\\
$T(0,0,4,1)\cong T(\lambda(\frac{-3}{4},\frac{-13}{8},\frac{3}{4}))$.
\end{tabular} 
\end{center}

\begin{figure}[h]
{\scalebox{0.55}{
\begin{tikzpicture}
\draw [step=1,thin, gray] (-5,-1) grid (5,2); 

\draw (0,0) circle [radius=0.3]; 
\filldraw (0,0) circle [radius=0.15]; 
\filldraw  [thick, fill=black] (0,1) circle [radius=0.17] ; \filldraw  [thick, fill=black] (1,0) circle [radius=0.17] ;
\filldraw  [thick, fill=black] (1,1) circle [radius=0.17] ; \filldraw  [thick, fill=black] (2,0) circle [radius=0.17] ;
\filldraw  [thick, fill=black] (2,1) circle [radius=0.17] ; \filldraw  [thick, fill=black] (3,0) circle [radius=0.17] ; 
\filldraw  [thick, fill=black] (3,1) circle [radius=0.17] ; \filldraw  [thick, fill=black] (4,1) circle [radius=0.17] ; 
\filldraw  [thick, fill=black] (-1,0) circle [radius=0.17] ; \filldraw  [thick, fill=black] (-1,1) circle [radius=0.17] ; 
\filldraw  [thick, fill=black] (-2,0) circle [radius=0.17] ; \filldraw  [thick, fill=black] (-2,1) circle [radius=0.17] ; 
\filldraw  [thick, fill=black] (-3,0) circle [radius=0.17] ; 

\filldraw  [thick, fill=white] (4,0) circle [radius=0.17] ; \filldraw  [thick, fill=white] (-4,0) circle [radius=0.17] ;
\end{tikzpicture}
}}
\caption{MCMs arising from consistent dimer models for Type 8b} 
\label{8b}
\end{figure}

As we mentioned, there are two consistent dimer models associated with $R$, and one of them is isoradial and the other is not isoradial. 
The following figures are an isoradial one and the associated quiver. (The fraction on each arrow is an $\sfR$-charge which indicates this dimer model is isoradial.)

%%%%%%%%%%%%%%%%%%%%%%%%%%%%%%%%%%%
%%%%%%%%%%%%Type_8b-1%%%%%%%%%%%%%%%%%%
%%%%%%%%%%%%%%%%%%%%%%%%%%%%%%%%%%%

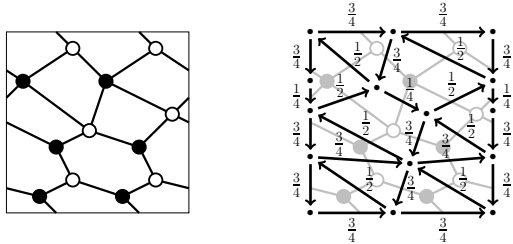
\begin{center}
\begin{tikzpicture} 
\node (DM) at (0,0) 
{\scalebox{0.4363}{
\begin{tikzpicture}
%vertex
\node (B1) at (1,0.5){$$}; \node (B2) at (3.5,0.5){$$}; \node (B3) at (1.5,2){$$};  \node (B4) at (4,2){$$};
\node (B5) at (0.5,4){$$}; \node (B6) at (3,4){$$};
\node (W1) at (2,1){$$}; \node (W2) at (4.5,1){$$}; \node (W3) at (2.5,2.5){$$};  \node (W4) at (5,3){$$}; 
\node (W5) at (2,5){$$}; \node (W6) at (4.5,5){$$};
\draw[thick]  (0,0) rectangle (5.5,5.5);

%edge
\draw[line width=0.068cm]  (B1)--(W1); \draw[line width=0.068cm]  (B1)--(0,0.75); \draw[line width=0.068cm]  (B1)--(1.5,0); 
\draw[line width=0.068cm]  (B2)--(W1); \draw[line width=0.068cm]  (B2)--(W2); \draw[line width=0.068cm]  (B2)--(4,0); \draw[line width=0.068cm]  (W1)--(B3);
\draw[line width=0.068cm]  (W2)--(B4); \draw[line width=0.068cm]  (W2)--(5.5,0.75); \draw[line width=0.068cm]  (B3)--(W3);
\draw[line width=0.068cm]  (B3)--(0,2.75); \draw[line width=0.068cm]  (W3)--(B4); \draw[line width=0.068cm]  (W3)--(B5); 
\draw[line width=0.068cm]  (W3)--(B6); \draw[line width=0.068cm]  (W4)--(B4); \draw[line width=0.068cm]  (W4)--(B6); 
\draw[line width=0.068cm]  (B5)--(W5); \draw[line width=0.068cm]  (B6)--(W5); \draw[line width=0.068cm]  (B6)--(W6); 
\draw[line width=0.068cm]  (W4)--(5.5,2.75); \draw[line width=0.068cm]  (W4)--(5.5,3.5); \draw[line width=0.068cm]  (B5)--(0,3.5);
\draw[line width=0.068cm]  (B5)--(0,4.33); \draw[line width=0.068cm]  (W5)--(1.5,5.5); 
\draw[line width=0.068cm]  (W6)--(4,5.5); \draw[line width=0.068cm]  (W6)--(5.5,4.33);
%black
\filldraw  [ultra thick, fill=black] (1,0.5) circle [radius=0.2] ;\filldraw  [ultra thick, fill=black] (3.5,0.5) circle [radius=0.2] ;
\filldraw  [ultra thick, fill=black] (1.5,2) circle [radius=0.2] ;\filldraw  [ultra thick, fill=black] (4,2) circle [radius=0.2] ;
\filldraw  [ultra thick, fill=black] (0.5,4) circle [radius=0.2] ;\filldraw  [ultra thick, fill=black] (3,4) circle [radius=0.2] ;
%white
\draw  [ultra thick,fill=white] (2,1) circle [radius=0.2] ;\draw  [ultra thick, fill=white] (4.5,1) circle [radius=0.2] ;
\draw  [ultra thick,fill=white] (2.5,2.5)circle [radius=0.2] ;\draw  [ultra thick, fill=white] (5,3) circle [radius=0.2] ;
\draw  [ultra thick,fill=white] (2,5)circle [radius=0.2] ;\draw  [ultra thick, fill=white] (4.5,5) circle [radius=0.2] ;
\end{tikzpicture}
}}  ;

\node (Qui) at (4,0) 
{\scalebox{0.4363}{ 
\begin{tikzpicture}
%vertex
\node (B1) at (1,0.5){$$}; \node (B2) at (3.5,0.5){$$}; \node (B3) at (1.5,2){$$};  \node (B4) at (4,2){$$};
\node (B5) at (0.5,4){$$}; \node (B6) at (3,4){$$};
\node (W1) at (2,1){$$}; \node (W2) at (4.5,1){$$}; \node (W3) at (2.5,2.5){$$};  \node (W4) at (5,3){$$}; 
\node (W5) at (2,5){$$}; \node (W6) at (4.5,5){$$};
\node (Q0) at (3,1.5){{\Large$\bullet$}}; \node (Q1a) at (0,3.1){{\Large$\bullet$}}; \node (Q1b) at (5.5,3.1){{\Large$\bullet$}}; 
\node (Q2) at (2,3.8){{\Large$\bullet$}}; \node (Q3) at (3.5,3){{\Large$\bullet$}}; \node (Q4a) at (0,4){{\Large$\bullet$}}; \node (Q4b) at (5.5,4){{\Large$\bullet$}}; 
\node (Q5a) at (2.5,0){{\Large$\bullet$}}; \node (Q5b) at (2.5,5.5){{\Large$\bullet$}}; \node (Q6a) at (0,0){{\Large$\bullet$}}; \node (Q6b) at (5.5,0){{\Large$\bullet$}}; 
\node (Q6c) at (5.5,5.5){{\Large$\bullet$}}; \node (Q6d) at (0,5.5){{\Large$\bullet$}}; \node (Q7a) at (0,1.7){{\Large$\bullet$}}; \node (Q7b) at (5.5,1.7){{\Large$\bullet$}}; 

%\node (Q0) at (3,1.5){{\Large$0$}}; \node (Q1a) at (0,3.1){{\Large$1$}}; \node (Q1b) at (5.5,3.1){{\Large$1$}}; 
%\node (Q2) at (2,3.8){{\Large$2$}}; \node (Q3) at (3.5,3){{\Large$3$}}; \node (Q4a) at (0,4){{\Large$4$}}; \node (Q4b) at (5.5,4){{\Large$4$}}; 
%\node (Q5a) at (2.5,0){{\Large$5$}}; \node (Q5b) at (2.5,5.5){{\Large$5$}}; \node (Q6a) at (0,0){{\Large$6$}}; \node (Q6b) at (5.5,0){{\Large$6$}}; 
%\node (Q6c) at (5.5,5.5){{\Large$6$}}; \node (Q6d) at (0,5.5){{\Large$6$}}; \node (Q7a) at (0,1.7){{\Large$7$}}; \node (Q7b) at (5.5,1.7){{\Large$7$}}; 
%\draw[thick]  (0,0) rectangle (5.5,5.5);

%edge
\draw[lightgray, line width=0.068cm]  (B1)--(W1); \draw[lightgray, line width=0.068cm]  (B1)--(0,0.75); \draw[lightgray, line width=0.068cm]  (B1)--(1.5,0); 
\draw[lightgray, line width=0.068cm]  (B2)--(W1); \draw[lightgray, line width=0.068cm]  (B2)--(W2); \draw[lightgray, line width=0.068cm]  (B2)--(4,0); 
\draw[lightgray, line width=0.068cm]  (W1)--(B3);
\draw[lightgray, line width=0.068cm]  (W2)--(B4); \draw[lightgray, line width=0.068cm]  (W2)--(5.5,0.75); \draw[lightgray, line width=0.068cm]  (B3)--(W3);
\draw[lightgray, line width=0.068cm]  (B3)--(0,2.75); \draw[lightgray, line width=0.068cm]  (W3)--(B4); \draw[lightgray, line width=0.068cm]  (W3)--(B5); 
\draw[lightgray, line width=0.068cm]  (W3)--(B6); \draw[lightgray, line width=0.068cm]  (W4)--(B4); \draw[lightgray, line width=0.068cm]  (W4)--(B6); 
\draw[lightgray, line width=0.068cm]  (B5)--(W5); \draw[lightgray, line width=0.068cm]  (B6)--(W5); \draw[lightgray, line width=0.068cm]  (B6)--(W6); 
\draw[lightgray, line width=0.068cm]  (W4)--(5.5,2.75); \draw[lightgray, line width=0.068cm]  (W4)--(5.5,3.5); \draw[lightgray, line width=0.068cm]  (B5)--(0,3.5);
\draw[lightgray, line width=0.068cm]  (B5)--(0,4.33); \draw[lightgray, line width=0.068cm]  (W5)--(1.5,5.5); 
\draw[lightgray, line width=0.068cm]  (W6)--(4,5.5); \draw[lightgray, line width=0.068cm]  (W6)--(5.5,4.33);
%black
\filldraw  [ultra thick, draw=lightgray, fill=lightgray] (1,0.5) circle [radius=0.2] ;\filldraw  [ultra thick, draw=lightgray, fill=lightgray] (3.5,0.5) circle [radius=0.2] ;
\filldraw  [ultra thick, draw=lightgray, fill=lightgray] (1.5,2) circle [radius=0.2] ;\filldraw  [ultra thick, draw=lightgray, fill=lightgray] (4,2) circle [radius=0.2] ;
\filldraw  [ultra thick, draw=lightgray, fill=lightgray] (0.5,4) circle [radius=0.2] ;\filldraw  [ultra thick, draw=lightgray, fill=lightgray] (3,4) circle [radius=0.2] ;
%white
\draw  [ultra thick,draw=lightgray,fill=white] (2,1) circle [radius=0.2] ;\draw  [ultra thick, draw=lightgray,fill=white] (4.5,1) circle [radius=0.2] ;
\draw  [ultra thick,draw=lightgray,fill=white] (2.5,2.5)circle [radius=0.2] ;\draw  [ultra thick, draw=lightgray,fill=white] (5,3) circle [radius=0.2] ;
\draw  [ultra thick,draw=lightgray,fill=white] (2,5)circle [radius=0.2] ;\draw  [ultra thick, draw=lightgray,fill=white] (4.5,5) circle [radius=0.2] ;

\draw[->, line width=0.08cm] (Q0)--node[midway,xshift=5pt,yshift=12pt] {{\huge$\frac{1}{2}$}}(Q1a); 
\draw[->, line width=0.08cm] (Q0)--node[midway,xshift=8pt,yshift=0pt] {{\huge$\frac{3}{4}$}}(Q5a); 
\draw[->, line width=0.08cm] (Q0)--node[midway,xshift=-5pt,yshift=12pt] {{\huge$\frac{3}{4}$}}(Q7b); 
\draw[->, line width=0.08cm] (Q1a)--node[midway,xshift=-3pt,yshift=12pt] {{\huge$\frac{1}{2}$}}(Q2); 
\draw[->, line width=0.08cm] (Q1a)--node[midway,xshift=-12pt,yshift=0pt] {{\huge$\frac{3}{4}$}}(Q7a); 
\draw[->, line width=0.08cm] (Q1b)--node[midway,xshift=12pt,yshift=0pt] {{\huge$\frac{3}{4}$}}(Q7b); 
\draw[->, line width=0.08cm] (Q2)--node[midway,xshift=7pt,yshift=10pt] {{\huge$\frac{1}{4}$}}(Q3); 
\draw[->, line width=0.08cm] (Q2)--node[midway,xshift=12pt,yshift=5pt] {{\huge$\frac{1}{2}$}}(Q6d); 
\draw[->, line width=0.08cm] (Q3)--node[midway,xshift=-8pt,yshift=5pt] {{\huge$\frac{3}{4}$}}(Q0); 
\draw[->, line width=0.08cm] (Q3)--node[midway,xshift=-7pt,yshift=12pt] {{\huge$\frac{1}{2}$}}(Q4b); 
\draw[->, line width=0.08cm] (Q4a)--node[midway,xshift=-12pt,yshift=0pt] {{\huge$\frac{1}{4}$}}(Q1a); 
\draw[->, line width=0.08cm] (Q4b)--node[midway,xshift=12pt,yshift=0pt] {{\huge$\frac{1}{4}$}}(Q1b); 
\draw[->, line width=0.08cm] (Q4b)--node[midway,xshift=15pt,yshift=7pt] {{\huge$\frac{1}{2}$}}(Q5b); 
\draw[->, line width=0.08cm] (Q5a)--node[midway,xshift=0pt,yshift=-15pt] {{\huge$\frac{3}{4}$}}(Q6b); 
\draw[->, line width=0.08cm] (Q5a)--node[midway,xshift=15pt,yshift=5pt] {{\huge$\frac{1}{2}$}}(Q7a); 
\draw[->, line width=0.08cm] (Q5b)--node[midway,xshift=10pt,yshift=0pt] {{\huge$\frac{3}{4}$}}(Q2);
\draw[->, line width=0.08cm] (Q5b)--node[midway,xshift=0pt,yshift=15pt] {{\huge$\frac{3}{4}$}}(Q6c); 
\draw[->, line width=0.08cm] (Q6a)--node[midway,xshift=0pt,yshift=-15pt] {{\huge$\frac{3}{4}$}}(Q5a); 
\draw[->, line width=0.08cm] (Q6b)--node[midway,xshift=10pt,yshift=8pt] {{\huge$\frac{1}{2}$}}(Q0);
\draw[->, line width=0.08cm] (Q6c)--node[midway,xshift=12pt,yshift=0pt] {{\huge$\frac{3}{4}$}}(Q4b); 
\draw[->, line width=0.08cm] (Q6d)--node[midway,xshift=-12pt,yshift=0pt] {{\huge$\frac{3}{4}$}}(Q4a); 
\draw[->, line width=0.08cm] (Q6d)--node[midway,xshift=0pt,yshift=15pt] {{\huge$\frac{3}{4}$}}(Q5b); 
\draw[->, line width=0.08cm] (Q7a)--node[midway,xshift=-12pt,yshift=0pt] {{\huge$\frac{3}{4}$}}(Q6a); 
\draw[->, line width=0.08cm] (Q7a)--node[midway,xshift=-18pt,yshift=14pt] {{\huge$\frac{3}{4}$}}(Q0); 
\draw[->, line width=0.08cm] (Q7b)--node[midway,xshift=12pt,yshift=0pt] {{\huge$\frac{3}{4}$}}(Q6b);
\draw[->, line width=0.08cm] (Q7b)--node[midway,xshift=10pt,yshift=8pt] {{\huge$\frac{1}{2}$}}(Q3);
\end{tikzpicture}
}} ;
\end{tikzpicture}
\end{center}
Again, by the results in \cite[subsection~5.12]{Nak}, we can see rank one MCM modules arising from an isoradial dimer model is actually conic, 
and non-conic divisorial ideals arise from a consistent dimer model which is not isoradial.

%%%%%%%%%%%%%%%%%%%%%%%%%%%%%%%%%%%%%%%%%%%%%%%%%%%%%%%%%%%%%%%%%%

\subsection*{Acknowledgements}
The author would like to thank Professor Mitsuyasu Hashimoto for valuable discussions about the Frobenius push-forward of toric singularities, 
Professor Kevin Tucker for informing the reference \cite{Bru}, Professor Hailong Dao for helpful comments about rank one MCM modules, 
and Kazunori Matsuda for stimulating discussions about toric singularities. 

The author is supported by Grant-in-Aid for JSPS Fellows No. 26-422.

%%%%%---reference---%%%%%

\end{document}